\documentclass[12pt,a4paper]{article}
\usepackage{euscript,fullpage,amsmath,amssymb}
\usepackage[T1]{fontenc}
\usepackage{enumerate}
\usepackage{tikz}
\usepackage{graphicx}
\usepackage{float}
\usetikzlibrary{patterns}
\usepackage{setspace}
\usepackage{amssymb}
\usepackage{subcaption}
\usepackage[utf8]{inputenc}

\newtheorem{theo}{Theorem}

\newtheorem{proposition}{Proposition}
\newtheorem{definition}{Définition}
\newtheorem{remark}{Remark}

\def\D{d_{\mu}^f}
\def\pn{\tilde{P}_n}
\title{Extreme value distributions of observation recurrences}
\date{}

\begin{document}

 \maketitle

\begin{center}
\begin{center}
Th. Caby \footnote{Aix Marseille Univ, CNRS, Centrale Marseille, I2M, Marseille, France E-mail: caby.theo@gmail.com},
D. Faranda\footnote{Laboratoire des Sciences du Climat et de l'Environnement, UMR 8212 CEA-CNRS-UVSQ,
IPSL and Universit\'e Paris-Saclay, 91191 Gif-sur-Yvette, France and
London Mathematical Laboratory, 8 Margravine Gardens, London, W6 8RH, UK. Email: davide.faranda@lsce.ipsl.fr.},
S.\ Vaienti\footnote{Aix Marseille Universit\'e, Universit\'e de Toulon, CNRS, CPT, 13009 Marseille, France. E-mail: vaienti@cpt.univ-mrs.fr.},
P. Yiou\footnote{Laboratoire des Sciences du Climat et de l'Environnement, UMR 8212 CEA-CNRS-UVSQ,
IPSL and Universit\'e Paris-Saclay, 91191 Gif-sur-Yvette, France. E-mail: pascal.yiou@lsce.ipsl.fr.}

\end{center}

 \begin{abstract}
We study analytically and numerically the extreme value distribution of observables defined along the temporal evolution of a dynamical system. The convergence to the Gumbel law of observable recurrences gives information on the fractal structure of the image of the invariant measure by the observable. We provide illustrations on idealized and physical systems.
\end{abstract}

\end{center}
\tableofcontents
\section {Introduction}

\subsection{A general overview} Extreme value theory (EVT) has been 
used in dynamical systems in the last years to quantify the probability of visiting a small set in the phase space, which constitutes a {\em rare event}. With this approach, the asymptotic statistics of hitting times and of the number of visits \cite{ei} in  small sets can be described. Methods based on EVT and more generally on the recurrence properties of chaotic systems  have found applications in climate science \cite{nature, messori, FF,nature2}. Quantifying the recurrence properties of weather patterns via dynamical indicators has proven useful to solve a  number of issues in climate and atmospheric sciences. In \cite{nature} the recurrence properties of the North-Atlantic sea-level pressure fields have been studied. A  number of instantaneous metrics that track rarity, predictability and persistence of atmospheric jet states and circulation patterns have been derived starting from quantities defined in the framework of EVT for dynamical systems, e.g. the local dimensions and the extremal index. In ~\cite{brunetti,rodrigues} the same metrics have been used to classify and evaluate the dynamical consistence of state-of-the art climate models in representing the atmospheric dynamics.   The impact of climate change of the atmospheric dynamical features was identified through shifts of the local dimensions between 1850 and 2100, in various datasets (observations, ensembles of scenario climate model simulation) \cite{nature2}. A critical discussion of the methods used in these studies is available in \cite{dq,ei}. To justify them, one needs to work with data sampled from the original high dimensional system, while experimentalists often have  access to a lower dimensional representation of the underlying attractor through measurements. A first approach to recover information on the underlying system from observations is to use embedding techniques, which is allowed by Takens' theorem \cite{taken}. Thanks to the theory of extreme value distribution applied to observables  developed  in this paper, we are able to propose an alternative technique and we will propose an application to atmospheric sciences. On a more general ground,  the aim of our work is to study the statistics of recurrences of smooth observables in chaotic dynamical systems. We will state some general results that could be applied in a wide range of situations. Our basic inspirations were the works of \cite{JR, hunt, MB}, where the authors developed different theoretical ideas and tools to  derive, among others,  recurrence rates for  observations and  compute them for various dynamical systems.

\subsection{Salient results of the paper}

 \begin{enumerate}

 \item {\bf Section 2} puts the basis of EVT for observations. We look at the distribution of the maximum of a sequence of random variables obtained by evaluating a vector valued observable  along the orbit of a dynamical system  and approaching a limiting value of the observable itself (the target set). We 
 obtain rigorously  a limit distribution of Gumbel type by using a perturbation theory applied to dynamical systems of hyperbolic type. 
 
 \item  An extremal index (EI)  modulates the  limit distribution, by adding a factor to the Gumbel law. This EI is related to the frequency of the occurrences (visits to the target sets), which is  interpreted as a clustering of the orbits. The EI becomes smaller than one when the target set exhibits periodic patterns. In {\bf Section 3} we  first provide general formulas for the EI for a large class of one-dimensional expanding maps and non-invertible observables. Then we  show that the observable could generate several coexisting clusters and we  explicitly compute the EI in a few cases. 
 
 \item 
 The numerical approach to the limit distribution via the Generalized Extreme Value (GEV) distribution, allows us to estimate the local properties of the image  measure.
 {\bf Section 4} is partially devoted to a brief exposition of the Hunt and Kaloshin theory of prevalent spaces in relation with the point-wise dimension of image measures. We therefore study in details two examples, the baker map and the product of two Cantor sets, and show that a few quite simple observables are not prevalent. This means that the dimension of the image measure is not integer (that of the ambient space), but smaller or larger and coinciding with that of the underlying attractor for the dynamics. The theoretical results were supported by numerical computations using the EVT techniques. This, combined with a suitable choice of the observable, is therefore a  very efficient tool to describe the fine geometric structure of the limit sets of the dynamics. 

\item We go beyond the Gumbel law in {\bf Section 5}, by  studying the statistics of the number of visits of the observable in the neighborhood of a value of interest. This is the point process associated to the distribution of the first hitting time, and we  show that it is either purely Poisson distributed or it deviates from the usual   Poly\`a-Aeppli distribution, which characterizes the point process when the rare set is around a periodic point. A particular example is studied in detail and a limit compound Poisson  distribution is exhibited {\em via} its generating  function and a recursive formula for the probability mass function. Application to  climate data  shows a  compound Poisson distribution, despite the relative modest length   of the time series  and the unavoidable approximations in their detection.
 
 \item In {\bf Section 6} we consider  what happens when the dynamical system and the observable are randomly perturbed. We show with analytical and numerical arguments, that if the perturbation of the map produces a smooth stationary measure or the observable changes randomly but staying prevalent, then the dimension of the image measure  becomes integer. Stability behaviors are also discussed.
 
 \item We then move to open systems in {\bf Section 7}  by considering in the phase space the presence of absorbing regions ({\em holes}),  where the orbits could be trapped and disappear forever. Nevertheless and under general conditions, a fractal repeller  survives and it is possible to study the recurrence properties of observables defined in a neighborhood of such a repeller.
 
 \item {\bf Section 8} gives a geometrical interpretation to our results and shows that our approach can be used to compute the hitting time statistics in the neighborhood of hypersurfaces embedded in the phase space of the system. Applications to fractal sets are also given.
 
 \item The experimental and numerical computation of the local dimension by EVT shows a discrete variability of such dimensions, even if they are constant almost everywhere (at least when they exist almost surely with respect to  ergodic measures).  The presence of those (large) deviations, is revealed by the non-linearity of the so-called spectrum of generalized dimensions (the {\em free-energy function} of the process), which are accessible to analytic and numerical computations. In {\bf Section 9}  we  treat  the large deviations of the dimensions of the image measure and discuss how those deviations are influenced by the choice of observable.  
 
 \item We quoted in section 1.1 the embedding technique as a tool to reconstruct the    attractor by considering the iterates of a unidimensional projection of the dynamics. When considering enough delay coordinates, the dimension of the attractor becomes accessible. In {\bf Section 10} we  propose an alternative approach that  allows us to have access to the dimension of the attractor by using directly observational data. In particular, this is possible when the dimensionality of the observations is larger than the information dimension of the underlying system. To achieve this either we  dispose of a vector-valued observable, or we could use a scalar observable to construct several images just by composing with the dynamics.  In some sense, the delay coordinate observable used in embedding techniques is a particular case of the smooth observables that we  consider. 
\end{enumerate}

\section{The formal approach}
We now introduce the basic concepts on Extreme Value Theory and apply them to  a sequence of observations. The stationary random process that  arises is then studied with a perturbative spectral technique, which allows us  to prove directly  the convergence to the Gumbel law.

\subsection{Basics on Extreme Value Theory}

Let us consider a dynamical system $(X, T, \mu)$, where $T$ acts on the measurable space $X$ and preserves the invariant probability measure $\mu.$ In the following we will consider $X$ as a compact subset of $\mathbb{R}^n$, ($n\ge 1$) and we put the Borel $\sigma$-algebra on it. We take $f: X\rightarrow \mathbb{R}^l$ a measurable function, called the {\em observable}; it will play a fundamental role in this paper, and additional assumptions on its regularity will be progressively added. 

Let us now construct the new measurable function
\begin{equation}\label{OB}
\phi(x)=-\log(\text{dist}(f(x)-f(z))),
\end{equation}
where $z$ is given in $X.$ This function has values in $\mathbb{R}\cup\{+\infty\}$ and achieves a global maximum at the pre-images of $f(z),$ where it is precisely infinite. With $\text{dist}$ we take a distance defining the metric on $\mathbb{R}^l.$ Consider the maximum of the process $\{\phi \circ T^k\}_{k>0}$, namely
\begin{equation}\label{MAX}
M_n(x)=\max\{\phi(x), \dots, \phi(T^{n-1}(x)\}
\end{equation}

and the distribution
\begin{equation}\label{PRO}
\mu(M_n\le u_n)=\mu(\phi\le u_n,\dots, \phi\circ T^{n-1}\le u_n).
\end{equation}

\begin{definition}
We say that we have an extreme value law for $M_n$ if there is a non-degenerate distribution function $H: \mathbb{R}^{+}\rightarrow [0,1]$ and for every $\tau>0$ there exists a sequence of levels $u_n=u_n(\tau), n\in \mathbb{N},$ such that 
\begin{equation}\label{SCA}
n\mu(\phi>u_n)\rightarrow \tau,  \ \text{as} \ n\rightarrow \infty.
\end{equation}
and for which the following holds:
$$
\mu(M_n\le u_n)\rightarrow H(\tau), \ \text{as} \ n\rightarrow \infty.
$$
\end{definition}

\begin{remark}
We name Eq. (\ref{SCA}) the {\bf Assumption F}: it allows us to avoid a degenerate limit for the distribution of $M_n.$ We will see later on  that the perturbative spectral  technique prescribes {\bf Assumption F} in a very natural way.
\end{remark}

Notice that Eq. (\ref{SCA}) is equivalent to
\begin{equation}\label{SCA2}
n \mu(x\in X, f(x)\in B(f(z), e^{-u_n}))\rightarrow  \tau,
\end{equation}
where $B(a, r)$ denotes the ball of radius $r$ centered at the point $a$ in the metric given by the chosen distance.

By introducing the {\em image measure} $f_*\mu$ defined as

\begin{equation}\label{IM}
f_*\mu(A)=\mu(f^{-1}A),
\end{equation}

where $A$ is any Borel set in $\mathbb{R}^l$, we can equivalently rewrite Eq. (\ref{SCA2}) as

\begin{equation}\label{SCA3}
n\ f_*\mu(B_{n,z})\rightarrow \tau,
\end{equation}

where we set

\begin{equation}
B_{n,z}:=B(f(z), e^{-u_n}) \,\mathrm{ and } \ C_{n,z}:=B(f(z), e^{-u_n})^c.
\end{equation}

The superscript $A^c$ is the complementary set of $A$ in $X.$

\begin{remark}

The presence of the observable imposes some natural conditions on the combined choice of $f$ and $T$ if we want to satisfy Eq. (\ref{SCA2}). For instance if $f$ is locally constant in the neighborhood of the target point $z$ and $\mu$ is not atomic in $z,$ we see immediately that Eq. (\ref{SCA2}) cannot hold for large $n.$ A less trivial example is given by the direct product map $T$ on the unit square defined by

$$T(x,y)=
\begin{cases}
 2x, \ x\in[0,1/2]; 1-2x, \ x\in [1/2, 1],\\
 ay, \ 0<a<1, y \in [0,1].
\end{cases}$$

 This map preserves the product of the Lebesgue measure on the $x$-axis times the Dirac mass at $0$ on the $y$-axis. If we now take the observable $f(x,y)=y$ and the target point in $(0,0),$ we see that the set $T^{-1}[B(f(0,0), e^{-u_n})]$ is a strip of length $1$ and of width  $e^{-u_n}$ on the square and the measure of this strip will be
 $1$ for any $n$.
 \end{remark}
Notice that if the observable  $f$ is not locally constant in the neighborhood of the target point  and the image measure is not atomic we can always choose a sequence $u_n$ verifying for each $n$: $n\ f_*\mu(B_{n,z})=\tau$. We will see in the next section, in particular the scaling (\ref{SCA4}), that $u_n$ is an affine function of the variable $y:=-\log \tau$ which can be written as: 
\begin{equation}\label{UU}
    u_n=\frac{y}{a_n}+b_n, \ a_n>0.
\end{equation}
When the sequence $\mu(M_n\le u_n)=\mu(a_n(M_n-b_n)\le y)$
converges to a non-degenerate distribution function $G(y)$, in the point of continuity of the latter, then we have an extreme value law. The starting point of EVT,  related to the affine choice for
the sequence $u_n$, is that such a $G(y)$ could be only of three types, called Gumbel, Fr\'echet and
Weibull (see \cite{led} for a general account of the theory). One of the main goal of this paper is to show that for the particular observable Eq. (\ref{OB}), we will get the Gumbel law, see  {\bf Proposition \ref{pro1}}. The scaling (\ref{SCA4}) shows that the parameters $a_n$ and $b_n$ are expressed in terms of the local dimension of the image measure, in Eq. (\ref{LLL}). It would be therefore useful to have access to those parameters. In this regard, we begin to notice that 
 the distribution function of the form $\mu(M_n\le y)$ is modeled for $n$ sufficiently large, by the so-called {\em generalized extreme value (GEV)} distribution \cite{PI}, which is a function depending upon three parameters $\xi\in \mathbb{R},$ (the tail index), $\kappa\in \mathbb{R},$ (the location parameter) and  $\hat{\sigma}>0,$ (the scale parameter):
$$
\text{GEV}(y; \kappa, \hat{\sigma}, \xi)=\exp\left\{-\left[1+\xi\left(\frac{y-\kappa}{\hat{\sigma}}\right)\right]^{-1/\xi}\right\}.
$$
The location parameter $\kappa$ and the scale parameter $\hat{\sigma}$ are scaling constants in place of $b_n$ and $a_n$. The idea is now to use a block-maxima approach (see Section 4.1) and fit our unnormalised data to a GEV distribution; for that it will be necessary to find a linkage among $a_n, b_n, \kappa$ and $\hat{\sigma}.$ For observables $\phi$ producing the Gumbel law, it has been shown in \cite{GEV}, that for $n$ large we have
\begin{equation}\label{rtt}
a_n\sim \frac{1}{\hat{\sigma}}; \ b_n\sim \kappa; 
\end{equation}
moreover the shape parameter $\xi$ tends to zero. The systematic use of this approach  from Section 4, will allow us to compute the local dimensions of the image measures and give therefore a numerical and experimental support to the theoretical results: this is another relevant aspect of our work. We finish this section by giving another definition.
\begin{definition}
 We say that the process $\{\phi\circ T^k\}_{k>0},$ for the observable (\ref{OB}),   has an Extremal Index (EI) $0 \le \theta\le 1,$  if we have a Gumbel distribution as
 $$
 \mu(M_n\le u_n)\rightarrow e^{-\theta \tau}, \ n\rightarrow \infty, 
 $$
with the sequence $u_n$ verifying {\bf Assumption F.}

\end{definition}
As we anticipated above, Gumbel's law is the limiting distribution for the maxima. The next sections will be devoted to the analytic computation of the EI. Besides rigorous estimates, we will also proceed to numerical computations. The EI is less than one when clusters of successive recurrences happen, which is the case, for instance, when the target point $z$ is periodic.   In our paper \cite{ei} we showed that the usual algorithms to compute  the  EI    have strong limitations   when  clusters  of higher order are present, and a new technique was proposed which consists in computing the first five $q_k$ terms in the expansion of $\theta$, see formula (\ref{EI}). We used this technique for the numerical estimates of the EI all along the paper.

\subsection{The perturbative spectral approach}
In order to apply the aforementioned perturbative spectral technique, we suppose that the system $(X, T, \mu)$ is REPFO ({\em Rare events Perron-Frobenius operators}) according to the terminology introduced by G. Keller \cite{GK, KL}. The definition of a REPFO dynamical system is quite technical even if its assumptions are verified in several situations when the system is uniformly hyperbolic or expanding. We must detail those assumptions because they  impose new constraints on the choice of the observable $f$. The basic object is the transfer  (Perron-Fr\"obenius)  operator $P$ associated to the map $T$. This operator acts on a suitable Banach space $(\mathcal{B}, ||\cdot||),$ equipped with a second (weak) norm $|\cdot|$ for which the closed unit ball of $(\mathcal{B}, ||\cdot||)$ is $|\cdot|$-compact. 

The Banach space is a space of functions or of distributions. We will mostly  treat non-invertible maps and in this case $\mathcal{B}$ will be the space of bounded variation ($BV$) functions and the weak norm will be the $L^1$ norm with respect to the Lebesgue measure. We will also consider invertible maps and, in this case, $\mathcal{B}$ is a space of distribution and we defer to \cite{demers} for a nice presentation  of those spaces or to  \cite{vv} for an easy description in the context of EVT. To make the exposition simpler we will suppose that $\mathcal{B}$ is the space of $BV$ functions and the weak norm is the space of integrable functions with respect to the Lebesgue measure $\text{Leb}$.\footnote{Sometimes, especially in the integral, we will write $dx$ instead of $d\text{Leb}$.} 

We will see below that the operator $P$ is slightly perturbed to get a sequence of operators $\pn$ which converge to $P$ in a sense that we are going to precise: for the moment we retain that $\pn$ is defined as $\pn (g)=P({\bf 1}_{W_n}g), \ g \ \text{and} \ {\bf 1}_{W_n}g\in \mathcal{B}$, where the Lebesgue measure of the measurable set $W_n$ goes to one when $n\rightarrow \infty,$ (see below for the explicit construction of such a $W_n;$ its complementary set $W_n^c$ should be interpreted as a {\em hole} with vanishing measure, not necessarily with vanishing diameter).  The following four items  define precisely what a REPFO system is: they are taken from \cite{GK} and slightly modified to our situations:

\begin{itemize}

\item {\bf A1} The unperturbed operator $P$ is quasi-compact:  this means, in particular,  that  $1$ is a simple isolated eigenvalue and there are no other eigenvalues on the unit circle. This implies the existence of a  unique mixing invariant measure $\mu$ for $T$ which is absolutely continuous with respect to $\text{Leb}$ with the density $h.$ 

    \item {\bf A2} There are constants $\beta, D>0,$ such that $\forall n$ sufficiently large, $\forall g\in \mathcal{B}$ and $\forall k\in \mathbb{N}$ we have ({\em Lasota-Yorke inequalities}):

\begin{eqnarray}
|\pn^k g|\le D |g|, \\
||\pn^k g||\le D \beta^k||g||+D|g|.
\end{eqnarray}

\item {\bf A3}  We can bound the weak norm of  $(P-\pn)g,$ with $g\in \mathcal {B},$ in terms of the norm of $g$ as:

$$|(P-\pn)g|\le \chi_n ||g||,$$

where $\chi_n$ is a monotone upper semi-continuous sequence converging to zero; this is called the {\em triple norm} estimate.

\item {\bf A4} If we put for $g\in \mathcal{B}$

\begin{equation}
\eta_{n}:=\sup_{||g||\le 1}|\int P(g\ {\bf 1}_{W_n^c})d\text{Leb}|,
\end{equation}

we must show that

\begin{equation}\label{R1}
\lim_{n\rightarrow \infty}\eta_n=0,
\end{equation}

\begin{equation}\label{R2}
\eta_{n}||P({\bf 1}_{W_n^c}h)||\le \text{const} \ \mu(W_n^c).
\end{equation}

\end{itemize}
 We now associate to the space of $BV$ functions a
 uniformly expanding  endomorphism $T$ of the unit interval and  preserving the absolutely continuous invariant mixing measure $\mu$ with density $h$.
The transfer operator has now a simple definition; for $v\in L^1(\text{Leb})$ and $w\in L^{\infty}(\text{Leb}),$ we have
$$
\begin{cases}
\int P(v) w dx= \int v w\circ T dx ,\\
 P(h)=h.
 \end{cases}
$$
Using this duality relation, the distribution in Eq. (\ref{MAX}) reads
\begin{equation}\label{K1}
\mu(M_n\le u_n)=\int ({\bf 1}_{C_{n,z}}\circ f)(x) \cdots ({\bf 1}_{C_{n,z}}\circ f)(T^{n-1}x) h(x) dx=
\int (\tilde{P}_n^nh)(x) dx.
\end{equation}
where
\begin{equation}\label{O1}
\tilde{P}_n g:= P({\bf 1}_{C_{n,z}}\circ f \ g), \ g\in \text{BV}.
\end{equation}

In the case of hyperbolic diffeomorphisms, we have a slightly different formula, since the operator   acts on measures, not on functions.
When $n\rightarrow \infty,$ the preceding assumptions
allow us to express the largest eigenvalue of $\tilde{P}_n,$ say $\chi_n,$ in terms of the largest eigenvalue of the unperturbed operator, which is $1$, as: 
$\chi_n = 1-(\theta\Delta_n+o(\Delta_n)),$ where $\Delta_n =
\mu(f^{-1}(B_{n,z})).$  The quantity $\theta$ is formally defined as $\theta=1-\sum_{k=0}^{\infty}q_k,$ and the $q_k$ are given by the following limits, {\em when they exist}:
\begin{equation}\label{sca}
q_k=\lim_{n\rightarrow \infty}q_{k,n}, \ \text{where} \ q_{k,n}:=\frac{\int (P-\tilde{P}_n)\tilde{P}^k_n(P-\tilde{P}_n)(h)dx}{\Delta_n}.
\end{equation}
The operator  $\tilde{P}_n$ now decomposes as the sum of a projection along the one dimensional
eigenspace associated to the eigenvalue $\chi_n$  and an operator with a spectral radius exponentially
decreasing to zero and which can be neglected in the limit of large $n.$\footnote{This is precisely what quasi-compactness means.}   Remembering
this, writing $ \tilde{P}_nh_n=\chi_nh_n,$ with $h_n$ converging to $h$   in the 
$L^1(\text{Leb})$  norm and replacing
into the right hand side of Eq. (\ref{K1}) and after a few manipulations we get, by neglecting
higher order terms:
$$
\mu(M_n\le u_n)\approx e^{-\theta n\Delta_n}.
$$
The product  $n\Delta_n=nf_{*}\mu(B_{n,z})$  is now controlled by {\bf Assumption F}, which allows us to get a limiting distribution. The justification of the previous statements is a direct application
of Keller's theory \cite{GK} which gives the following

 \begin{proposition}\label{pro1}
 Let us suppose that $(X, T, \mu)$ is  a REPFO system. Suppose  moreover that
 the {\bf Assumption F}  holds. Then
\begin{equation}\label{EVT1}
\mu(M_n\le u_n)\underset{n\to\infty}\rightarrow e^{-\theta \tau},
\end{equation}
where the extremal index $\theta$ is defined as
\begin{equation}\label{EI}
\theta=1-\sum_{k=0}^{\infty}q_k,
\end{equation}
where
\begin{equation}\label{cucu}
q_k:=\lim_{n\rightarrow \infty} q_{k,n}
\end{equation}
The quantities $q_k$ are given by the limit (\ref{sca}) where the quantities $q_{k,n}$ are equivalently expressed as:
\begin{equation}\label{qk}
q_{k,n}= \frac{1}{f_*\mu(B_{n,z})}\mu(f^{-1}B_{n,z}\cap T^{-1}(f^{-1}B_{n,z})^c\cdots \cap T^{-k}(f^{-1}B_{n,z})^c \cap T^{-(k+1)}f^{-1}B_{n,z}).
\end{equation}
 \end{proposition}

 {\em Comments.} As we said above the proof of this proposition follows immediately from Keller's theory, see also  our previous works  \cite{D1,D2,dq,ei, vv}. Three issues deserve to be discussed. The first two deal with the possibility to give  examples which fit with the REPFO assumptions. Whenever $f$ is the identity function, the aforementioned references give a large class of examples. The problem now is the presence of the observable $f$ which could affect the hypothesis {\bf A2}-{\bf A4}. The third issue concerns the computation of the extremal index.
 \begin{enumerate}
     \item
  A particular attention must be drawn to the Lasota-Yorke inequalities {\bf A2}  which has to do  with the characteristic function of sets of the type $f^{-1}(C_{n,z}),$ which could have a geometric shape quite different from  balls. We should guarantee that the Banach norm of these sets is computable and allows to get the desired inequalities. This will be the case for all the systems with associated observables which we will be  treated analytically in this paper. We will in fact consider the observables  $f$ as  continuous and  local  $C^1$ functions and in this case the Lasota-Yorke inequalities for the  perturbed operators  can be proved using the arguments in \cite{BBVV}, Lemma 2.6, or \cite{LM}, Lemma 7.4. For the baker map (see section 4), we defer to the paper \cite{vv}, section 3.1.
 \item The second isuue concerns the Assumptions {\bf A3} and {\bf A4}.  They will follow if we could prove that the
 $L^1(\text{Leb})$ norm of  $(\tilde{P}_n-P)g$, with $g\in BV,$ is  bounded by the $BV$ norm of $g$ and the image-Lebesgue ($f_*\text{Leb}$) measure of $B_{n,z}.$ \\
 We have
 \begin{eqnarray}
 \int |(\tilde{P}_n-P)g|h(x)dx & =  & \int |P({\bf 1}_{B_{n,z}}(f(x))g(x))|h(x)dx \nonumber\\
  & \le & ||h||_{\infty}\int P({\bf 1}_{B_{n,z}}(f(x))|g(x)|)dx \nonumber\\
 & \le & ||h||_{\infty}\int {\bf 1}_{B_{n,z}}(f(x))|g(x)|dx \nonumber\\
 & \le & ||h||_{\infty}||g||_{BV} f_*\text{Leb}(B_{n,z}),
 \end{eqnarray}

since both $h$ and $g$ are in $BV$ and the infinity norm is bounded by the $BV$ norm $||\cdot||_{BV}.$ The perturbative theorem requires finally that $f_*\text{Leb}(B_{n,z})\le \text{constant} \ f_*\mu(B_{n,z}),$ which is surely true if the density $h$ is bounded from below: we will tacitly assume it if necessary.
\item Finally we should check the existence of the limits (\ref{cucu}) to give the fundamental  expression of Eq. (\ref{EI}) for the extremal index. Note that  the Poincar\'e recurrence theorem  implies that $\sum_{k=0}^{\infty}q_{k,n}=1$; therefore whenever $q_k$ exists,  the extremal index is at most $1$. The quantities $q_{k,n}$ have a simple geometrical interpretation: they give the conditional measure of the points that are at the beginning in the set  $f^{-1}B_{n,z},$ are iterated outside it for the next $k$ times, and finally return to it at the $k+1$ iteration. As we will argue below, in particular in section 4, this structure of the $q_{k,n}$ allows us to compute them explicitly in several situations, or guess their possible behavior.\\
 \end{enumerate}

We now define the local dimensions of the image measure. We put
\begin{eqnarray}\label{LLL}
\underline{d}_{\mu}^f(x):=\liminf_{r\rightarrow 0}\frac{\log f_*\mu (B(f(x),r))}{\log r},\\
\overline{d}_{\mu}^f(x):=\limsup_{r\rightarrow 0}\frac{\log f_*\mu (B(f(x),r))}{\log r}.
\end{eqnarray}

  Whenever $\underline{d}_{\mu}^f(x)=\overline{d}_{\mu}^f(x)=d_{\mu}^f(x),$ we will say that the image measure $f_*\mu$ is {\em exact dimensional.}\\

  {\bf Notations:} Sometimes instead of $\D(x)$ we will use the notation $\D(f_0)$, meaning that the pointwise dimension is computed in the point $f_0=f(z)$ without specifying the value of $z.$ When the measure is exact dimensional we will simply write $\D$ as the almost sure value. We will also use the symbol $\D$ to denote what we presume to be the almost sure value of the image measure in a few numerical computations for which the measure $\mu$ could only be reconstructed numerically. This especially concerns the last chapter.\\

 We now suppose that $d_{\mu}^f(x)$ exists and express  Eq. (\ref{SCA3}) as the scaling:
\begin{equation}\label{SCA4}
f_*\mu(B_{n,z})\sim e^{-u_n\ d_{\mu}^f(z)}\sim \tau/n.
\end{equation}

The result on the Gumbel law given by Eq. (\ref{EVT1}) could be reformulated by: if $f$ is an observable on the space $X$, and $f(z)$ is the value at a given point $z$, then {\em the probability that the distance between $f(T^nx)$ and $f(z)$ after $n$ iterations is less than $(\frac{\tau}{n})^{\D}$ for the first time, is approximately $e^{-\theta \tau}.$}\\

Let us come back to the results in \cite{JR}. The equality between the recurrence rate for the observable defined as
\begin{equation}\label{rf}
R_f(z)=\lim_{r\to0} \frac{\log \inf\{k\in \mathbb{N}^*: f(T^kz)\in B(f(z),r)\}}{\log r},
\end{equation}

and $d_{\mu}^f(z)$ is proven for a class of systems with superpolynomial decay of correlations and such that the image measure $f_*\mu$ is exact dimensional.
 In the spectral theory, the property of superpolynomial decay of correlation is  strengthened   by the presence of the spectral gap for the transfer operator, which implies exponential decay of correlations. We point out that our approach is slightly different from the one of \cite{JR}, in the sense that we get a recurrence rate for hitting times instead of return times and its distribution for shrinking target sets. For this reason we will not further elaborate on the connections with the quantity in Eq. (\ref{rf}).




\section{The extremal index}

The extremal index $\theta$ is usually considered as a measure of {\em clustering}, whenever several and repeated occurrences take place in the ball $B_{n,z}$. For the usual observable $-\log \text{dist}(x,z)$, this happens around periodic points for the map $T$. When it comes to recurrence of observables, some clustering can also occur when $z$ is a periodic point. We now show that the extremal index for observables reveals new interesting features.

We start with a simple example.

Take a real observable $f$ defined on the unit interval such that in any point where it is defined, the derivative is bounded below away from zero and above from infinity.

Let us first consider the case on an {\em invertible} $f$
and take $T$ as a uniformly expanding map of the interval which has $z=f^{-1}(f_0),$ where $f_0=f(z),$ as a fixed point and is continuous in such a point together with the density of the absolutely continuous invariant measure $h$. Then we have

\begin{equation}\label{eq:q0n}
q_{0,n}= \frac{1}{f_*\mu(B_{n,z})}\mu(f^{-1}B_{n,z} \cap T^{-1}f^{-1}B_{n,z}).
\end{equation}

At this point we can repeat the standard argument (see, for example, section 4.2 of \cite{targets}) to get immediately that

$$\theta= 1- \frac{1}{|T'(z)|}.$$

We now take a {\em non-invertible} $f$. In particular we suppose $f$ has two branches: $f_1, f_2.$\\
Suppose the ball $B_{n,z}$ is again centered at a point $f_0=f(z)$ and the point $z_{1}:=f^{-1}_{1}(z)$ is the inverse point of $f$ such that
$$Tz_{1}=z_1.$$
Moreover suppose that the other pre-image $z:=f^{-1}_{2}(z)$ is not periodic for $T$.

In Eq. (\ref{eq:q0n}) for the $q_{0,n}$ above, only the pre-images by $T$ of the set $f^{-1}_{1}(B_{n,z})$  matter in the computation of the EI, but we have to take into account the relative ratio of the measure of  $f^{-1}_1(B_{n,z}), f^{-1}_2(B_{n,z})$ in the denominator. These measures are obtained by pulling back the Lebesgue measure of $B_{n,z}$ with the reciprocal images of $f$, which amounts to multiply the length of $B_{n,z}$ with the reciprocal of the derivative of $f$ in the pre-images of $f_0$, and multiply what we get by the density $h$ in such pre-images.  In conclusion we have
\begin{equation}
\theta=1-\frac{1}{|T'(z_1)|}\frac{1}{1+\frac{h(z) |f'(z_1)|}{h(z_1) |f'(z)|}} .
\end{equation}\label{form1}
The preceding argument can be  generalized to give an exact formula for the $q_k$. As we will see, the existence of several pre-images of the ball $B_{n,z}$ could generate multiple clusters coexisting with different degrees of periodicity.

\begin{proposition}\label{PPP1}

Let us suppose that $T$ is a uniformly expanding map as above and the observable $f$ is differentiable with a derivative bounded away from zero and infinity. Fix $z$ in the unit interval $M$ and put $f_0=f(z);$ suppose also that $f$ is a finite-to-one map. Consider the set of the pre-images $w$ of $f_0,$ one of them being $z,$ and suppose that they do not belong  to the countable union of the pre-images of the boundary points of the domains of local injectivity of $T$ and that the invariant density $h$ is continuous in such points. Consider the set $$A_k=\{w\in M: f(w)=f_0, f(Tw)\neq f_0,\dots, f(T^kw)\neq f_0, f(T^{k+1}w)=f_0\}.$$

 When $A_k=\emptyset$, then $q_k=0$. Conversely, whenever $A_k$ is finite and non-empty, we have
\begin{equation}\label{gen}
q_k=\underset{w\in A_k}\sum \frac{1}{|T^{(k+1)}(w)'|}\frac{1}{1+ \frac{|f'(w)|}{h(w)}\underset{y\in B^{\omega}_k}\sum \frac{h(y)}{|f'(y)|}},
\end{equation}
where $B_k^{\omega}=\{y\in M: f(y)=f_0\}\setminus \{w\}$.\\

The extremal index is obtained by $$\theta=1-\sum_{k=0}^{\infty} q_k.$$

\end{proposition}

We point out that having fixed the center $z$ of the ball $B_{n,z}$ and having $f$ a finite number of pre-images, there are only finitely many points in $A_k$  and consequently finitely many terms in the sum $\sum_{k=0}^{\infty} q_k$. Moreover we could relax the global  assumption on $f$ by asking that $f$ be $C^1$ in $z$ and the pre-images of $z$ as the next example will require.

Let us give two examples. In the first consider the map $T(x)=3x$-mod $1$. Then take a point $a>\frac12$ which is not periodic for $T$ (these points yield a full Lebesgue measure), and consider a piece-wise continuous straight line with two branches $(f_1, f_2),$ $f_1$ passing through  the  points $(0,0)$ and $(a,1)$, and $f_2$ through the points $(a,1)$ and $(1,0).$ The equations are
$$\begin{cases}
f_1(x)= x/a,\\
f_2(x)= \frac{1}{a-1}\ x- \frac{1}{a-1}.
\end{cases}$$
We choose a point $z_1$ that is a fixed point of $T$ and $z_2$ that have the same image by $f$ but is not periodic by $T$. We take
$$\begin{cases}
z_1=\frac12,\\
z_2= f_2^{-1}(f(z_1))=(a-1)/2a+1.
\end{cases}$$
We can choose $a$ so that $z_2$ is irrational. In this case, $z_2$ is not periodic for $T$ and the trajectory starting from $z_2$ will not pass through $z_1$ which is rational. Therefore, we see easily from formula (\ref{gen}) that $q_k=0$ for $k>0$.
We have $h(z_1)=h(z_2)=1$ and
$$\begin{cases}
|f'(z_1)|=1/a,\\
|f'(z_2)|= \frac{1}{|a-1|}.
\end{cases}
$$
Therefore
$$\theta=1-q_0= 1- \frac{1}{3\left(1+\frac{|a-1|}{a}\right)}.$$

 We checked this formula numerically for various values of $a$. For example, taking $a=2/\pi$, we find a numerical value of $0.788$ against  a theoretical value of $0.7878$. We used the estimate $\hat\theta_5$ introduced in \cite{ei}, which consists in estimating the $q_k$ terms up to the order $5$ and subtracting them from $1$.

We notice that in Eq. (\ref{gen}), the extremal index depends explicitly on the density of the invariant measure, which was constant in the example above.

We give another example where  $h(z_1)$ and $h(z_2)$ are different. We take the Hemmer map defined in $[-1,1]$ by $T(x)=1-2\sqrt{|x|}$. Its density is $h(x)=\frac12(1-x),$  \cite{hemmer}. The point $z_1=3-2\sqrt2$ is a fixed point of the map. We choose the point $z_2=-1/2$, which is not periodic and we take $f$ piecewise linear with different slopes: $f(x)=x$ for $-1\le x \le 0$ and $f(x)=-2x+11/2-4\sqrt2$ otherwise, so that $f(z_1)=f(z_2)=-1/2$. Eq. (\ref{form1}) gives
$$\theta=1-q_0=\frac{\sqrt{3-2\sqrt2}}{1+\frac3{4(\sqrt2-1)}}\approx 0.9104.$$

Our numerical computations confirm this result to the fourth digit with the estimate $\hat\theta_5$.\\

We have given a quite general formula for the one dimensional case, and it is apparent from it that the clustering structure can be quite complicated if the observable and the dynamics have some kind of compatibility. For this reason, giving a general formula for the extremal index in higher dimensional systems is out of the scope of this paper. We however believe that for large class of observables, no clustering is detected and the extremal index should be equal to $1$. This is confirmed by several numerical simulations that will be described in the next chapter.

\section{Phenomenology of the image measure }
\label{dddd}

We are now interested in estimating  the quantity $\D$ that appears in the distribution of maxima. This question has been partially answered by Rousseau and Saussol in \cite{JR}, in the case of smooth observables and measures $\mu$ that are absolutely continuous with respect to Lebesgue. In particular, Theorem 9 in  \cite{JR}  states that $d_{\mu}^f(z)$ exists $\mu$ almost everywhere, is integer valued and is equal to the rank of $Df(z)$ almost everywhere. For example, if $f$ has values in $\mathbb{R}$ and $\mu(\partial_xf(x) =0)=0$, $d_{\mu}^f(z)$ is equal to $1$ for $\mu$-almost any $z.$ It implies also that if $f$ is constant on some regions of the phase space of positive measure, $d_{\mu}^f(z)$ will be $0$ in that region. A first observation is that this result does not hold at some special points of the attractor. We now give an example where $d_{\mu}^f$ is not an integer.\\

 Consider the map $Tx=2x\mod1$ defined on the circle, $z=0$ and the observable $f(x)=x^a$, with $a>0$. Then we have: $$\mu(-\log|f(x)-f(0)|>u_n)=\mu(-\log|x^a|>u_n)=\mu(-\log|x|>\frac{u_n}{a})=\mu(B(0,e^{-\frac{u_n}{a}}))= 2e^{-\frac{u_n}{a}}.$$ 
 Therefore, $$d_{\mu}^f(0)=\lim_{n\to\infty} \frac{\log 2e^{-\frac{u_n}{a}} }{\log e^{-u_n}}=  1/a.$$ Depending on the value of $a$, this quantity can be non integer and either smaller or larger than 1.\\

In many physical applications, the measure is not smooth, but has a (multi)fractal structure. This happens for chaotic  dynamics in neuroscience and climate science \cite{brain,dq}. Being able to compute the value of $d_{\mu}^f$ in such situations is of crucial importance to describe the statistics of recurrence of the observable. The simplest case to consider is when the observable $f$ is a diffeomorphism from $\mathbb{R}^k$ to $\mathbb{R}^k$ (where $k$ is the dimension of the ambient space): the image of the invariant set by $f$ is then a deformation of the original attractor which preserves  its local structure. We therefore expect that $d_{\mu}^f(z)=\gamma_{\mu}(z)$, the pointwise dimension at the point $z$. Let us now remind the definition of these local dimensions, since they will be used later on.  Consider the limits
\begin{eqnarray}\label{LLDD}
\gamma_{\mu}^-(z)=\liminf_{r\rightarrow 0}\frac{\log \mu (B(z,r))}{\log r}\\
\gamma_{\mu}^+(z)=\limsup_{r\rightarrow 0}\frac{\log \mu (B(z,r))}{\log r}.
\end{eqnarray}
They are called respectively the lower and upper pointwise dimensions of $\mu$ at $z$. If $\gamma_{\mu}^-(z)=\gamma_{\mu}^+(z)$, the common value $\gamma_{\mu}(z)$ is called the pointwise dimension of $\mu$ at $z.$ We defer to our paper \cite{dq} for a discussion of these pointwise dimensions with the associated references. 

 Most observables used in practice are not diffeomorphisms. The most general result concerning the local dimension of image measures is due to Hunt and Kaloshin, in particular Theorem 4.1 in \cite{hunt}. Before stating their theorem, we must recall the important notion of {\em prevalence} used in the aforementioned paper, see also \cite{preval}, \cite{prevalence} and \cite{wiki}.  We consider a real topological vector space $V$ and a  Borel-measurable subset $S$  of $V$. $S$ is said to be {\em prevalent} if there exists a finite-dimensional subspace $P$ of $V$, called the {\em probe} set, such that for all $v\in V$ we have $v + p\in S$ for $\text{Leb}_P$-almost all $p\in P$, where $\text{Leb}_P$ denotes the $P$-dimensional Lebesgue measure on $P$. In the case of interest for us, $V$ is the space of $C^1$ functions $f:\mathbb{R}^n\rightarrow \mathbb{R}^m.$ The notion of prevalence could  be thought as the analogue of almost everywhere in infinite dimensional spaces. We give a few properties and examples of prevalence to point out its significance. All prevalent subsets $S$ of $V$ are dense in $V$. Then, if we declare that  {\em almost every}  means that the stated property holds for a prevalent subset of the space in question, we have:
 \begin{itemize}
 \item almost every continuous function from the interval $[0, 1]$ into $\mathbb{R}$ is nowhere differentiable. Here, $V$ is  the space of continuous functions on the unit interval with the supremum norm topology;
 \item take now $V=L^{1}(dx)$ the space of Lebesgue summable functions on the unit interval. Then
almost every function $f\in V$   has the property that $\int_0^1 f(x) dx\neq 0.$
\item if $A$ is a compact subset of $\mathbb{R}^n$ with Hausdorff dimension $d_H$, $m \ge d_H$, and $1 \le k \le\infty$, then, for almost every $C^k$ function $f:\mathbb{R}^n\rightarrow \mathbb{R}^m,$ $f(A)$ also has Hausdorff dimension $d_H.$
\end{itemize}
Other  examples will now be  stated in terms of dimension of measures. We summarize them in the following theorem:
\begin{theo}\label{HUU} (Hunt and Kaloshin \cite{hunt})
\begin{itemize}
\item Let $\mu$ be a Borel probability measure on $\mathbb{R}^n$ with compact support. For a
prevalent set of $C^1$ functions (also, for almost every linear transformation) $f: \mathbb{R}^n\rightarrow \mathbb{R}^m$,
$$
\underline{d}_{\mu}^f(x)=\min(m, \gamma_{\mu}^-(x))
$$
for almost every $x$ with respect to $\mu$. If in addition $\D(x)$ exists for almost every $x$, then for
almost every $f$ the pointwise dimension of $f$ at $f(x)$ exists and is given by
$$
\D(x)=\min(m, \gamma_{\mu}(x))
$$
for almost every $x.$
\item Let $\mu$  be a Borel probability measure on $\mathbb{R}^n$ with compact support. If the
pointwise dimension $\gamma_{\mu}(x)$ exists and does not exceed $m$ for almost every $x$ with respect
to $\mu$, then for a prevalent set of $C^1$ functions (also, for almost every linear transformation)
$f:\mathbb{R}^n \rightarrow \mathbb{R}^m$, the information dimension of $f$, $D_1(f(\mu))$ exists and is given by the information dimension $D_1(\mu)$ of $\mu$:
$$
D_1(f(\mu))=D_1(\mu).
$$
\end{itemize}
\end{theo}
The information dimension of a measure $\mu$ is defined as the following limit, when it exists
\footnote{Otherwise one should turn to the $\liminf$ and $\limsup.$
 We defer to \cite{hunt} and \cite{dq} for the details.}:
 $$
 D_1(\mu):=\lim_{r\rightarrow 0}\frac{\int \log\mu(B(x,r))d\mu(x)}{\log r}.
 $$
 We  note that with the assumptions of the second item of the preceding theorem, we have
 $$
 D_1(f(\mu)=\int \D(x) df_*\mu(x); \ D_1(\mu)=\int \gamma_{\mu}(x)d\mu(x).
 $$
 An important class of measures are those called {\em exact dimensional}: they enjoy the property that
 $$
 \gamma_{\mu}(x)=D_1(\mu), \ x-\mu\ a.e.
 $$

 {\bf Notations.} We will call {\em typical} a point $x$ that belongs to the set of full measure giving $D_1(\mu).$
 Sometimes we will simply write $D_1$ instead of $D_1(\mu)$ if the measure $\mu$ is clear from the context; moreover and still for exact dimensional measures we will use the short $\mu(B(x,r))\approx r^d$ in  place of $\lim_{r\rightarrow 0}\frac{\log\mu(B(x,r))}{\log r}=d.$\\

  Several dynamical systems
with hyperbolic properties have an invariant measure that
is exact dimensional. It is enough that the limit defining the local dimensions exists almost everywhere and that the measure is ergodic to have exact dimensionality \cite{young}.  In these cases the information dimension can be expressed in terms of the Lyapunov exponents and of the metric
entropy.

\begin{remark}
In the rest of the section,  we will consider a few  cases where we  compute $\D$ and compare it with the conclusions of the Hunt and Kaloshin Theorem. We will see that non-prevalent observables arise very easily in simple examples. With  abuse of language we will say that an observable is  prevalent if it belongs to the prevalent space of the Hunt-Kaloshin Theorem. We declare that an observable is not prevalent whenever it does not satisfy the Theorem above and for almost all choices of the target point $x$ (typical points). Later on  (example of the product of two Cantor sets), we will show an example of observable that violates the Hunt-Kaloshin Theorem for a given point $x$. 
Even in that case we will say that the observable is not-prevalent.  
 \end{remark}
 
\subsection{The baker map}
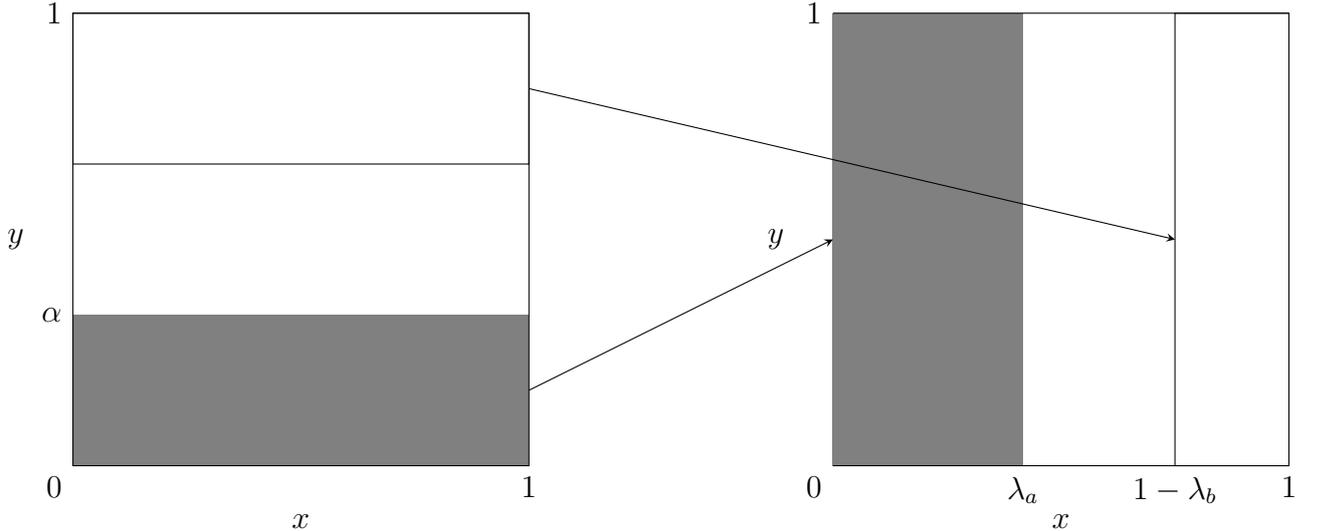
\begin{figure}
\begin{tikzpicture}
\fill [gray] (0,0) rectangle (6,2);
\draw (0,4) rectangle (6,6);
\draw (0,0) rectangle (6,6);
\fill [gray] (10,0) rectangle (12.5,6);
\draw (14.5,0) rectangle (16,6);
\draw (0,0) node[below left]{$0$};
\draw (0,6) node[left]{$1$} ;
\draw (-0.5,3) node[left]{$y$} ;
\draw (0,2) node[left]{$\alpha$} ;
\draw (3,-0.5) node[below]{$x$} ;
\draw (6,0) node[below]{$1$} ;
\draw (10,0) node[below left]{$0$} ;
\draw (16,0) node[below]{$1$} ;
\draw (10,6) node[left]{$1$} ;
\draw (9.5,3) node[left]{$y$} ;
\draw (13,-0.5) node[below]{$x$} ;
\draw (12.5,0) node[below]{$\lambda_a$} ;
\draw (14.5,0) node[below]{$1-\lambda_b$} ;
\draw[-] (10,0) -- (16,0);
\draw[-] (10,6) -- (16,6);
\draw[->,>=stealth] (6,5) -- (14.5,3);
\draw[->,>=stealth] (6,1) -- (10,3);

\end{tikzpicture}
\caption{Action of baker's map on the unit square. The lower part of the square is mapped in the left part and the upper part in the right part.}
\label{tik}
\end{figure}

We start with the two dimensional dynamics defined by the baker map, whose fractal SRB measure has been extensively studied, \cite{ott}, \cite{VVV}. It is defined on the unit square $Q=[0,1]\times[0,1]$  by the equations
\begin{equation}
x_{n+1}=\left\{
\begin{aligned}
\lambda_ax_n, \text{                  } y_n<\alpha,\\
(1-\lambda_{b})+\lambda_{b}x_n, y_n>\alpha,\\
\end{aligned}
\right.
\end{equation}
and
\begin{equation}
y_{n+1}=\left\{
      \begin{aligned}
\frac{y_n}{\alpha}, y_n<\alpha,\\
\frac{y_n-\alpha}{1-\alpha}, y_n>\alpha,\\
\end{aligned}
\right.
\end{equation}\label{baker}

where $\alpha \in(0,1/2]$ and $\lambda_a+\lambda_b \le 1$.  The action of the map on the unit square is shown in figure \ref{tik}. The SRB measure is exact dimensional, and its information dimension  is given by \cite{ott}:
\begin{equation}\label{bakerinf}
D_1=1+D_{1,s}, 
\end{equation}
with
$$
D_{1,s}:=\frac{\alpha\log (\alpha^{-1})+(1-\alpha)\log((1-\alpha)^{-1})}{\alpha\log(\lambda_a^{-1})+(1-\alpha)\log (\lambda_{b}^{-1})}.
$$

The spectral approach to EVT used in section 2, applies to  baker's map \cite{vv}.
Let us first consider  the {\em mean value} observable defined as
$$f(x,y)=\frac{x+y}2.$$

To compute numerically the quantity $\D$, we generate a trajectory of $M=10^8$ points starting from a point $x$ chosen at random on the square and compute at each iteration the value of $\phi_z(T^ix)=-\log |f_0-f(T^ix)|$ \footnote{The orbit of $x$ will approach quickly the attractor and it will give the the right statistical information by definition of SRB (physical) measure.}. We then compute the empirical distribution of the maximum taken by $\phi_z$ over blocks of size $n=5\cdot 10^4$. The scale parameter $\hat\sigma$ of the GEV distribution is computed with a maximum likelihood estimate, using the Matlab function {\em gevfit} \cite{gevfit}. An estimate for $\D$ is then given by $1/{\hat\sigma}$. The estimates of $\D$ are then averaged over $10$ different trajectories. The results are displayed in table \ref{table1} (the error is the standard deviation of the results over the $10$ trajectories). Although the measure has a fractal structure, we found, for different values of $f_0$ and  $\alpha$, estimates    for $\D$ that are very close to $1$, as expected from the result of Hunt and Kaloshin. Since the proof of their theorem does not allow a clear geometrical understanding of what happens, we provide now an illustration and a heuristic explanation for that result.\\

Let us take a typical point $z=(z_1, z_2)$, not lying on the border of the square, such that  $f(z)=c$, $0<c<1$ and let $\varepsilon >0$. The points verifying  $|f(x,y)-c|\le \varepsilon,$ are those on the straight lines $\frac{x+y}{2}=s,\ c-\varepsilon<s<c+\varepsilon$. This defines a strip where each couple $(x,y)$ will meet infinitely many  vertical unstable leaves foliating the attractor. Then the ball $B(c, \varepsilon)$ is completely filled and it could be reasonable to argue that $\D=1.$ This would be true if the measure $f_*\mu$ is absolutely continuous, as prescribed in \cite{JR}. But there is no reason that $f_*\mu$ has such a property, if $\mu$ is not absolutely continuous. As Fig. 2 shows, we really found $\D=1$, which fits with Theorem \ref{HUU} and  suggests  that $f$ is prevalent. Before  giving  a rigorous direct proof of this fact, we point out  that the above example can be easily modified with a drastic change in the dimension of the image measure, which therefore exhibits a non-prevalent observable. We defer to the end of this section for such an example. 

  We therefore consider the  observable $f(x,y)=\frac{x+y}2$ and  the strip $\Sigma'_{\varepsilon}$ defined through $\frac{x+y}{2}=s,\ c-\varepsilon<s<c+\varepsilon$. We begin to remind that the SRB measure $\mu$ disintegrates along the vertical unstable leaves with absolutely continuous conditional measures (actually Lebesgue measures normalized to $1$), and with singular measures along the horizontal stable leaves, we will use them later on.  The unstable leaves $W_{u, \iota}$ are indexed by $\iota$ and counted by the counting measure $\zeta';$ as we said, the conditional measure $\mu_{u,\iota}$ is the linear Lebesgue measure of mass $1$. Then the SRB measure of the strip $\Sigma'_{\varepsilon}$ reads:
  $$
  \mu(\Sigma'_{\varepsilon})= \int \mu_{u,\iota}(\Sigma_{\varepsilon}'\cap W_{u,\iota})d\zeta'(\iota).
  $$
But $\mu_{u,\iota}(\Sigma_{\varepsilon}'\cap W_{u,\iota})=2\varepsilon$ and what are left in the integral above are an ensemble of unstable leaves of finite $\zeta'$ measure due to the affine term $c$ cutting the $y$-axis.  In conclusion $\mu(\Sigma'_{\varepsilon})\approx \varepsilon$ in agreement with
theorem \ref{HUU}.  We notice that the previous proof adapts easily to all affine observables of type $f(x,y)=ax+by+c$, provided that $b$ is different from zero.

 As an   interesting example of violation of prevalence, we  take $f$ as a multivariate Gaussian function maximized at the typical  point $z=(x_0,y_0)$, with a covariance matrix equal to the identity:
 \begin{equation}\label{MM}
f(x,y)=\frac1{2\pi}\exp\left(-\frac12\left((x-x_0)^2+(y-y_0)^2\right)\right).
\end{equation}
The set of points on $Q$ for which $|f(x,y)-f(x_0, y_0)|\le \varepsilon$ are the points belonging to the ball $B(z, 2\sqrt{\pi \varepsilon}):$
$$
0\le (x-x_0)^2+(y-y_0)^2\le 4\pi \varepsilon.
$$
Since the point $(x_0, y_0)$ is typical, we have $\D=\frac{D_1}{2}.$

\begin{figure}[h!]
    \centering
    \begin{subfigure}[t]{0.5\textwidth}
        \centering
\includegraphics[height=2.5in]{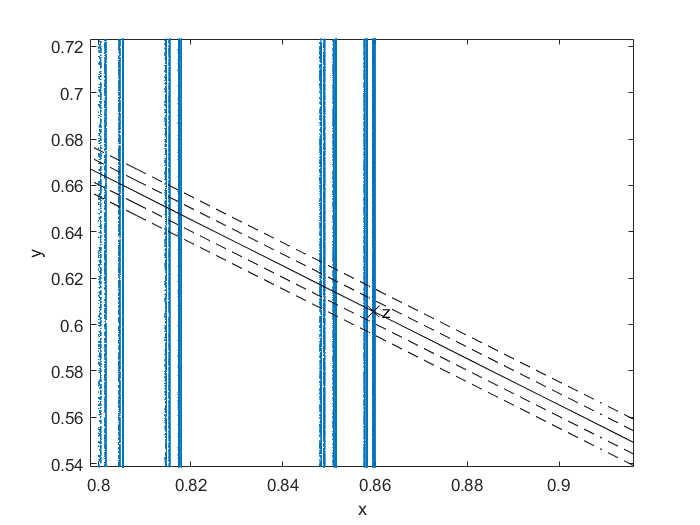}
\label{theta}
    \end{subfigure}%
    ~
    \begin{subfigure}[t]{0.5\textwidth}
        \centering
     \includegraphics[height=2.5in]{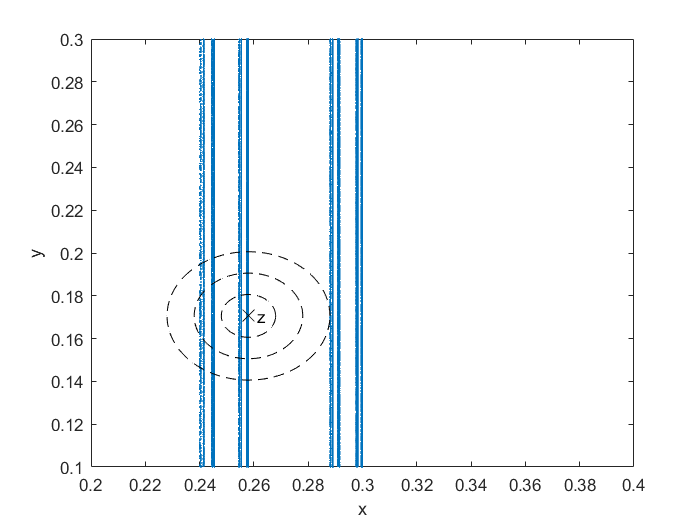}
     \label{d}
    \end{subfigure}
   \caption{Pictorial representation of the situation described in the main text for the observable $\frac{x+y}{2}$ (left) and a Gaussian centered at $z$ (right), in different regions of $Q.$ The baker attractor is depicted in blue, and the graphs $\{(x,y): f(x,y)=f(z) \pm \varepsilon\}$ are dotted lines. In both situations, these manifolds intersect the attractor an infinite number of times.}
   \label{fig}
\end{figure}

\begin{center}
\begin{tabular}{|c|c|c|c|}

  \hline
  &  $\alpha=1/5$  & $\alpha=1/4$ & $\alpha=1/3 $ \\
  \hline
  $f_0=0.1$ & $1.00 \pm 0.01$ & $1.00 \pm 0.02 $ &  $1.00 \pm 0.01$\\
  $f_0=0.3$ & $1.00 \pm 0.01$ & $1.00 \pm 0.01$ & $1.00 \pm 0.01$\\
  $f_0=0.8$ & $1.00 \pm 0.01$ & $1.00 \pm 0.01$ & $1.00 \pm 0.02$\\
\hline
\end{tabular}

\captionof{table}{Values of $\D(f_0)$ computed for the mean value observable, for different values of $\alpha$ and $f_0$. We took $\lambda_a=0.3$, $\lambda_b=0.2$.}
 \label{table1}
\end{center}

As mentioned earlier, we expect to detect no clustering of high values for such generic observables and non-periodic $z$. The extremal index is computed using the estimate $\hat\theta_5$, using as a threshold the $0.999$-quantile of the observable distribution. Results are averaged over $10$ trajectories and are presented in table \ref{table2}.

\begin{center}
\begin{tabular}{|c|c|c|c|}

  \hline
 &  $\alpha=1/5$  & $\alpha=1/4$ & $\alpha=1/3 $ \\
  \hline
  $f_0=0.1$ & $1 \pm 0$ & $1 \pm 0$ & $1 \pm 0$ \\
  $f_0=0.3$ & $1 \pm 0$ & $1 \pm 0$ & $1 \pm 0$\\
  $f_0=0.8$ &$1 \pm 0$ & $1 \pm 0$ & $1 \pm 0$\\
  \hline

\end{tabular}

\captionof{table}{Values of $\theta$ computed for the mean value observable, for different values of $\alpha$ and $f_0$. We took $\lambda_a=0.3$, $\lambda_b=0.2$. The error of $0$ is the standard deviation of the estimates.}
 \label{table2}
\end{center}

To get  the quantity $\D$  different from $1$ on a set of full measure, we should take an observable with range at least in $\mathbb {R}^2$. 

We performed  numerical computations using the baker map with parameters $\alpha=1/3$, $\lambda_a=1/3$, $\lambda_b=1/4$ and by taking the observable $f(x,y)=(x,x^2+y^2)$.  For different points $z$ not lying on the $x-$axis, we find indeed  a value of $\D(z)$ that is close to the information dimension $=1.2682$, as it is computed from Eq. (\ref{bakerinf}). For the point $z=(0.9581,0.0612)$ for example, we find a local dimension equal to $1.26 \pm 0.03$. We used again the parameters $M=10^{8}$ and $n=2\cdot10^5$.\\

  As promised above, we  now give another example of a non-prevalent observable.  Let us take the function $f(x,y)=x$ and $f_0=f(z_1,z_2)=z_1=c$, where $(z_1, z_2)$ is a typical point.  We need to compute the scaling of the SRB measure of the vertical strip $\Sigma_{\varepsilon}:=\{(x,y)\in Q; |x-c|\le \varepsilon\}$. To this end, we disintegrate the SRB measure $\mu$   along the horizontal stable leaves. These measures can be seen as generated by an Iterated Function System (IFS) with two scales $\lambda_a, \lambda_b$ and two weights $\alpha, 1-\alpha,$ \cite{VVV}. We now evaluate the SRB measure of the strip $\Sigma_{\varepsilon}$, as:
\begin{equation}\label{rt}
 \mu(\Sigma_{\varepsilon})=\int \mu_{s,\nu}(\Sigma_{\varepsilon}\cap W_{s,\nu})d\zeta(\nu),
\end{equation}
 where $\mu_{s,\nu}$ is the conditional measure along the stable leaf $W_{s,\nu}$, indexed by $\nu$ and counted by the counting measure $\zeta.$ These conditional measures are the same on each $W_{s,\nu}$ and for almost all choices of $z_1$  they behave as exact dimensional fractal measure with the exponent given by the term $D_{1,s}$ in Eq. (\ref{bakerinf}): 
 $$\mu_{s,\nu}(\Sigma\cap W_{s,\nu})\approx \varepsilon^{D_{1,s}}.$$
 
  Since the counting measure $\zeta(Q)=1,$ we finally get $\mu(\Sigma_{\varepsilon})\approx \varepsilon^{D_{1,s}},$  which violates theorem \ref{HUU}\, because the exponent should be equal to $1.$

\subsection{The product of two Cantor sets}
As a second example which can be worked out analytically, we  consider  the cartesian product of two ternary Cantor sets $K\times K$ on the unit interval $I.$ The dynamics is generated by two independent iterated function system (see \cite{ott}), each of them defined by  two linear contractive maps $g_1, g_2$ with slope $1/3.$  On each factor $K$ we take a measurable
map, our dynamical system,  $T : K \rightarrow K,$ with $T(x)=g_i^{-1}(x),$ for
 $x\in g_i(K).$ The Cantor set  $K$ will be the invariant set for
the transformation $T.$ The invariant measure $\mu^{(2)}=\mu \times \mu$ is the product of the two invariant measures on the factor spaces. Each factor measure  is a balanced measure with two equal weights $1/2$, which means that for any Borel set $B$ on the unit interval we have $\mu(B)=\sum_{i=1}^2\frac12 \mu(g_i(B)).$
 All these measures are exact dimensional and the information dimension of $\mu^{(2)}$ is $D_1=\frac{2\log2}{\log3}\approx1.26.$\\
The spectral approach to EVT used in section 2, applies to this systems, see \cite{D2} and \cite{targets}.
As a first observable, we take the standard multivariate Gaussian function (\ref{MM}). If we take a typical point  $z:=(x_0, y_0),$ we can repeat the argument given above for the baker's map and found $\D=D_1/2,$ which shows that (\ref{MM}) is not prevalent. This result is confirmed by the numerical simulations for which we
 used the same algorithm as described earlier for the baker map, using the parameters $M=5\cdot10^8$ and $n=5\cdot10^4$. We used more data and a different size of blocks because stable estimates are difficult to get. We detected no clustering, as for the baker's map. The results for the estimates of $\D$ are  shown in table \ref{cantortable}. We found for different points $z$ a value for $\D$ close to $0.61$, which is comparable with  $D_1/2\approx 0.63.$  \\

 We now take $f(x,y)=x$  and look at the strip $\Sigma_{\varepsilon}=\{|x-c|\le \varepsilon\},$ where $c$ is chosen on a typical line $(c,y).$ Instead of disintegrating, we can now use Fubini's theorem since $\mu^{(2)}$ is a product measure. We have
 $$
 \mu^{(2)}(\Sigma_{\varepsilon})=\int \mu_x(\Sigma_{\varepsilon}\cap K)d\mu_y
 $$
 where we write $\mu_x$ (resp. $\mu_y$) for the factor measure on the Cantor set $K$ on the $x$-axis (resp. on the $y$-axis). The sectional measure $\mu_x(\Sigma_{\varepsilon}\cap K)$ is independent  of $y$ and since it is exact dimensional it yields  $\mu_x(\Sigma_{\varepsilon}\cap K)\approx \varepsilon ^{D_1/2},$ which finally gives $\D=D_1/2,$ showing that the observable is not prevalent. 
 The same results holds for the observable $f(x,y)=y.$ \\

 A less trivial observable that violates prevalence is given by $f(x,y)=y-x$ and $f_0=f(0,0)=0$. We warn the reader that  we are not sure that the point $(0,0)$ is typical, but we discuss this situation since it could arise in  concrete applications, in the same way the periodic points are negligible in measure but play an important role in recurrence.  We therefore have   to compute the logarithm of the  $\mu^{(2)}$ measure of  the $\varepsilon$ neighborhood $\Sigma_{\varepsilon}$ of the diagonal $\{x=y\}$ and compare it with the logarithm of $\varepsilon$. It is easy to check that we can restrict ourselves  to countable sequences  like $\varepsilon_n=\beta^n$, with $\beta>1$ and $n\rightarrow \infty$. In particular we now take  $\varepsilon_n= 3 \ 3^{-n}$. Then we have as above
 $
 \mu^{(2)}(\Sigma_{\varepsilon})=\int \mu_x(\Sigma_{\varepsilon}\cap K)d\mu_y.
 $\\
Each time $x\in K$ on the $x$-axis, $\Sigma_{\varepsilon}$ will meet the Cantor set $K\times K$ in the point $(x,y).$ We therefore evaluate the $\mu_y$ measure of the section $\Sigma_{\varepsilon}\cap K$ by splitting it over the $2^n$ cylinders of the $n$-th generation in the construction of the Cantor set along the $y$-axis. There will be at least one of these cylinders of $\mu_y$-measure $2^{-n}$ inside that section. Then
$$
\mu^{(2)}(\Sigma_{\varepsilon_n})=\int \mu_x(\Sigma_{\varepsilon_n}\cap K)d\mu_y\ge 2^{-n}.\footnote{Actually we should remove from it a small contribution $o(2^{-n})$ due to integration on the corners $(0,0), (1,1),$ which will not affect the final result.}
$$
With the prescribed choice of $\varepsilon_n$ we finally get
$$
\limsup_{n\rightarrow \infty}\frac{\log\mu^{(2)}(\Sigma_{\varepsilon_n})}{\log \varepsilon_n}\le \frac{\log 2}{\log3},
$$
which immediately implies $\overline{d}^f_{\mu^{(2)}}(0,0)\le \frac{\log 2}{\log 3}.$\\

\begin{figure}[h!]
    \centering
    \begin{subfigure}[t]{0.5\textwidth}
        \centering
\includegraphics[height=2.5in]{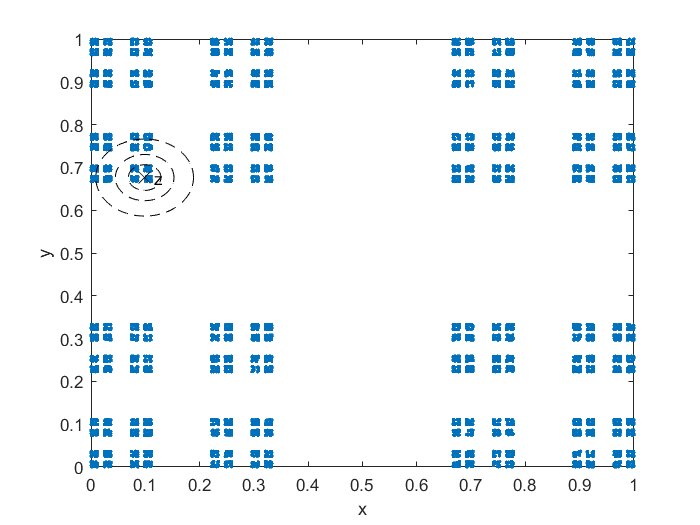}
\label{theta}
    \end{subfigure}%
    ~
    \begin{subfigure}[t]{0.5\textwidth}
        \centering
     \includegraphics[height=2.5in]{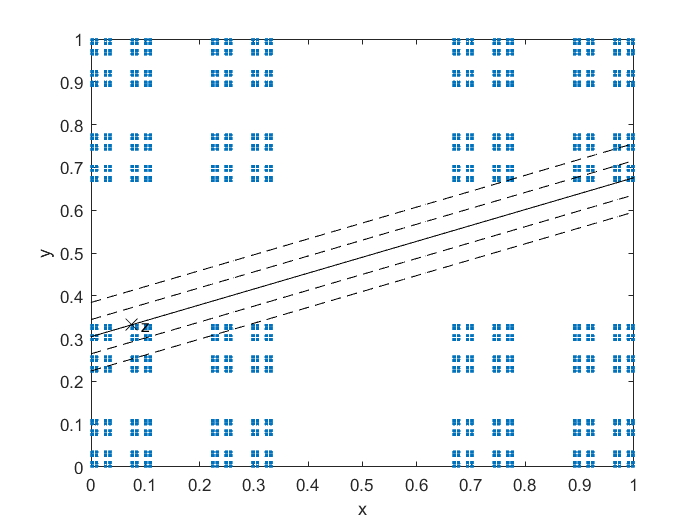}
     \label{f}
    \end{subfigure}
   \caption{For the product of Cantor sets, pictorial representation of the situation described in the main text for the gaussian observable (left) and a generic linear observable (right). The attractor is depicted in blue, and the graphs $\{(x,y): f(x,y)=f(z) \pm \varepsilon\}$ are dotted lines.}
   \label{K}
\end{figure}

 It is  difficult to prove directly that $\D=1$ as prescribed by the Hunt-Kaloshin theorem for the {\em almost sure prevalence affine functions} $f(x,y)=ax+by+c$ (see figure \ref{K} for a pictorial representation). Numerical experiments confirm such a behavior, although stable estimates are difficult to get. We chose the parameters $a,b$ and $c$ at random in the unit interval. For the point $z=(0.893,0.307)$, $a=0.557$, $b=0.6596$, $c=0.0046$, we find that $\D =0.993 \pm 0.02$. We used the same parameters as for the Gaussian observable.\\

\begin{center}
\begin{tabular}{|c|c|c|c|}
\hline
 $z$ &  $(0.994,0.0029)$  & $(0.6679,0.9914)$ & $(0.0861,0.2565)$\\
  \hline
$\D$ & $0.61 \pm 0.002$ & $0.60 \pm 0.002$ & $0.62 \pm 0.002$\\
  \hline

\end{tabular}\\

\captionof{table}{Values of $\D$ computed for the Gaussian observable, for different points $z$. The error is the standard deviation of the results.}
 \label{cantortable}
\end{center}

\subsection{The Lorenz system}
Let us now turn to a higher dimensional situation and consider the Lorenz 1963 system \cite{lorenz} that we reconstruct with the Euler method with step $h=0.01$. With this iterative procedure, the system can be seen as a discrete mapping, for which the developed theory is applicable. We chose an observable with image in $\mathbb{R}^5$. We tested several observables but the results are displayed for the observable $f(x,y,z)=(x^2+y^2,z,y+z,\pi yz,1/x)$. We find that the values of $\D(x,y,z)$ are all close to $D_1=2.06$. For the point $(-1.7323,8.9400,32.6818)$ for example, we have that $\D (x,y,z)=2.05 \pm 0.02$ using the parameters $M=10^8$ and $n=2\cdot 10^5$ (the results are averaged over 20 trajectories, and the error is the standard deviation of the results).\\

Instead, if we take a scalar observable, we find values very close to 1,
indicating  again that when the observable decreases the dimensionality, the fractal structure of the attractor is smoothed in the image measure (we are {\em supposing} here that the invariant measure is exact dimensional). These numerical results are in perfect agreement with the Hunt-Kaloshin Theorem.

\subsection{Conclusions} We conclude this section by pointing out the few examples which we found and do not verify the conclusions of the Hunt-Kaloshin Theorem. It happens when $f: \mathbb{R}^n \rightarrow \mathbb{R}^m$ with $m\le n.$ If the dimension of the attractor in $\mathbb{R}^n$ is larger than $m$, one expects to find $\D=m$ for a prevalent observable.  We exhibit several examples where, in the same circumstances, $\D<m.$ This shows that we are in presence of a non prevalent observable.

Another example of observable that does not belong to the prevalent set of the measure is a function whose Jacobian does not have a full rank on a set of positive measure, for absolutely continuous measures. This is a consequence of  theorem 9 in \cite{JR}. We emphasize that the image measure can have counter-intuitive properties. For instance, Rousseau \cite{thesejerome} gives the example of an image measure that is non atomic and yet $\D$ is $0$ on a set of positive measure. This example is built upon a Cantor set and the $C^{\infty}$ observable is defined as the limit of an iterative process.

\section{Statistics of visits for the observable}
\label{viss}
It can be interesting from a physical point of view to study the number of visits of the observable $f$ near a certain value $f_0=f(z)$. This problem is well understood in the framework of EVT. Let us consider the following counting function:
\begin{equation}\label{PP1}
N_n(t)=\sum_{l=1}^{\lfloor{\frac{t}{f*_\mu(B(f_0,r_n))}}\rfloor}{\bf 1}_{B(f_0,r_n)}(f(T^lx)),
\end{equation}
where the radius $r_n$ goes to $0$ when $n$ tends to infinity.
We are interested in the distribution
\begin{equation}\label{PP2}
\mu( N_n(t)=k), \ k\in \mathbb{N}
\end{equation}
when $n\rightarrow \infty.$
It has been proved (see for instance \cite{HV, FFT,FFT13}) that for $f=Id$ and when $z$ is not a periodic point, $\mu( N_n(t)=k)$ converge to the Poisson distribution $\frac{t^k e^{-t}}{k!}$, while for a periodic point of minimal period $p$, $\mu( N_n(t)=k)$  converges to the Poly\`a-Aeppli distribution, which is a particular kind of compound Poisson distribution.
Before continuing, we remind that  a probability measure $\tilde\nu$ on $\mathbb{N}_0$ is compound
Poisson distributed with parameters $t\lambda_\ell$, $\ell=1,2,\dots$,
if its generating function $\varphi_{\tilde\nu}$ is given by
$\varphi_{\tilde\nu}(z)=\exp\int_0^\infty(z^x-1)\,d\rho(x)$,
where $\rho$ is the measure on
 $\mathbb{N}$ defined by $\rho =\sum_\ell t\lambda_\ell\delta_\ell$,
with $\delta_\ell$  being the point mass at $\ell$.

 An important non-trivial compound Poisson distribution is the P\'olya-Aeppli
 distribution which holds when  the random variables given by the hitting times of the ball $B(f_0, r_n)$ is geometrically distributed, which implies that  $\lambda_\ell=(1-p)p^{\ell-1}$ for $\ell=1,2,\dots$, for some $p\in(0,1)$.
 In this case
 \begin{equation}\label{crty}
 \mathbb{P}(N_n(t)=k)\rightarrow e^{-pt}
\sum_{j=1}^kp^{k-j}(1-p)^j\frac{(pt)^j}{j!}\binom{k-1}{j-1}, \ n\rightarrow \infty,
\end{equation}
where $p$ is the extremal index. In particular $\mathbb{P}(W=0)=e^{-t}$.
In the case of $p=1$ this reverts  to the usual Poisson distribution.
For more general target sets, the limit law of $N_n(t)$ is given by a compound Poisson distribution when the extremal index is different from $1$, and by a pure Poisson distribution if no clustering occurs,  \cite{HV}, \cite{HV2}. We refer also to our paper \cite{ei} for a discussion of this matter and related references.

We now show that in presence of non-invertible observables $f$, we  get compound Poisson distributions which are not P\`olya-Aeppli.
\begin{proposition}
With the assumptions of Proposition \ref{PPP1}, suppose the ball $B(f_0, r_n),$ $f_0=f(z)$, has two  pre-images $B_{1,n}, B_{2,n}$, the first containing the periodic point $w_1=z$ of period $p_1,$ the second the periodic point $w_2$ of period $p_2.$
Then the distribution $N_n(t)$ is compound Poisson, but not P\`olya-Aeppli.\\

Proof: We notice that Eq. (\ref{PP1}) can be rewritten as
$$
N_n(t)=\sum_{l=1}^{\lfloor{\frac{t}{f*_\mu(B(z,r_n))}}\rfloor}{\bf 1}_{f^{-1}(B(f_0,r_n))}(T^lx)
$$
We can therefore apply the theory recently developed by \cite{HV2}, where entry times are considered for sets whose measure goes to zero. In our case those sets are the pre-images of the ball $B(f_0, r_n)$ and they are located around the points $w_i\in A_k, k\in \mathbb{N}$, where the set $A_k$ has been defined in Proposition 2; actually there are now only two pre-images.

 If we now refer to the theory in \cite{HV2} and in particular to Section 8.3 therein, we can easily compute the quantity $\tilde{\alpha}_l=\sum_i b_i^l,$ where $b_i^l:=\lim_{n\rightarrow \infty} \mu_{B_n}(\tau^{l-1}_{B_n}=i),$ being $B_n=B_{1,n}\cup B_{2,n}$ and $\tau_{B_n}^{l-1}$ is the $l-1$th return time into the set $B_n;$ with $\mu_A$ we intend the conditional measure to the set $A.$ For a given $l,$ only the terms $(l-1)b_1$ and $(l-1)b_2$  count  in the  sum defining $\tilde{\alpha}_l.$
By repeating the computation in Lemma 4 in \cite{HV2} we have
$$
\tilde{\alpha}_l=b_1^{l-1}\lim_{n\rightarrow \infty}\left[\frac{\mu(B_{1,n})}{\mu(B_n)}\right]+b_2^{l-1}\lim_{n\rightarrow \infty}\left[\frac{\mu(B_{2,n})}{\mu(B_n)}\right],
$$
where
$$
b_1=\lim_{n\rightarrow \infty} \left[\frac{\mu(B_{1,n}\cap T^{-p_1}B_{1,n})}{\mu(B_{1,n})}\right], \ b_2=\lim_{n\rightarrow \infty}\left[\frac{\mu(B_{2,n}\cap T^{-p_2}B_{2,n})}{\mu(B_{2,n})}\right]
$$
Notice that in the particular case we are considering and by repeating the computation in section 3 we have
\begin{eqnarray}
b_1=\frac{1}{|(T^{p_1})'(w_1)|}, \
b_2=\frac{1}{|(T^{p_2})'(w_2)|},\\
\mu_1:=\lim_{n\rightarrow \infty}\left[\frac{\mu(B_{1,n})}{\mu(B_n)}\right]=\frac{1}{1+\frac{h(w_2) |f'(w_1)|}{h(w_1) |f'(w_2)|}}, \\
\mu_2:=\lim_{n\rightarrow \infty}\left[\frac{\mu(B_{2,n})}{\mu(B_n)}\right]=\frac{1}{1+\frac{h(w_1) |f'(w_2)|}{h(w_2) |f'(w_1)|}}.
\end{eqnarray}
Moreover by recalling the definition of the quantities $q_k$ introduced in section 3 we have
$$
q_{p_1}=b_1\ \mu_1; \ q_{p_2}=b_2\ \mu_2.
$$
According to the theory developed in \cite{HV2}, the parameter $\lambda_l$ which we introduced before Eq. (\ref{crty}) to define the compound Poisson distribution is given by
 $$
 \lambda_l=\frac{\alpha_l-\alpha_{l+1}}{\alpha_1}, \
 \text{where} \
 \alpha_k=\tilde{\alpha}_k-\tilde{\alpha}_{k+1},
 $$
and $\alpha_1$ is the extremal index defined as the {\em reciprocal of the expected length of the clusters}:
$$
\sum_{k=0}^{\infty}k\lambda_k=\frac{1}{\alpha_1}.
$$
In our case and using the expression for the quantities introduced above we have:
$$
\alpha_1=1-(q_{p_1}+q_{p_2}).
$$
The latter  is an alternative way to define the extremal index,
which in the present situation  is  consistent with the formula found in section 3 for the extremal index $\theta$ using the spectral technique. We defer to our article \cite{ei} for a critical discussion of these  equivalent definitions. Moreover, putting $1-b_1\mu_1-b_2\mu_2=\theta,$ we have
$$
\lambda_l=\frac{b_1^{l-1}\mu_1(1-b_1)^2+b_2^{l-1}\mu_2(1-b_1)^2}{\theta}
$$
and for the generating function of the random variable given by the number of visits
\begin{eqnarray}\label{jj}
\phi(z) & = & \exp\left[\sum_{l=1}^{\infty}\theta t\lambda_l(z^l-1)\right] \nonumber\\
& = & e^{-\theta t}e^{t\mu_1(1-b_1)^2\frac{z}{1-zb_1}}
e^{t\mu_2(1-b_2)^2\frac{z}{1-zb_2}},
\end{eqnarray}
which gives a compound distribution different from the Poly\`a-Aeppli distribution.

\end{proposition}

We remind that deviations from the P\`olya-Aeppli distribution were exhibited in other situations, for instance  when the target set is a neighborhood of the diagonal in \cite{HV, ei} or  a neighborhood of periodic points where the map is not continuous in \cite{Ay}.\\

We now give a recursive formula that produces the distribution of $N_n(t)$. Let us denote $$a_1=t\mu_1(1-b_1)^2\  \text{and} \ a_2=t\mu_2(1-b_2)^2.$$
We first notice that
\begin{equation}\label{product}
\phi'(z)=\phi(z)\pi(z)',
\end{equation}
where
$$
\pi(z)=\frac{a_1z}{1-b_1z}+\frac{a_2z}{1-b_2z}.$$
We easily see that the $k$ derivatives of $\pi$ (for $k>0$) are given  by
\begin{equation}\label{deriv}
\pi^{(k)}(z)=\frac{k!a_1b_1^{k-1}}{(1-b_1z)^{k+1}}+\frac{k!a_2b_2^{k-1}}{(1-b_2z)^{k+1}}.
\end{equation}

Using the Leibniz formula for derivations, we have from Eq. (\ref{product}):
$$\phi^{(n)}(z)=\sum_{k=0}^{n-1}{n-1 \choose k}\phi^{(k)}(z)\pi(z)^{(n-k)}(z).$$

We now use this last formula and combine it with Eq. (\ref{deriv}) to obtain:
\begin{equation}\label{hhhh}
\phi^{(n)}(0)=\sum_{k=0}^{n-1}{n-1 \choose k}\phi^{(k)}(0)(n-k)!(a_1b_1^{n-k-1}+a_2b_2^{n-k-1}).
\end{equation}

Keeping in mind that from Eq. (\ref{jj}), $\phi(0)=e^{-\theta t},$ we can use  formula (\ref{hhhh}) to determine the probability that $N_n(t)=k$ by computing recursively the derivatives of the generating function $\phi$ at 0 and dividing by $k!$.\\

\begin{figure}[t]
\centering
\includegraphics[height=2.5in]{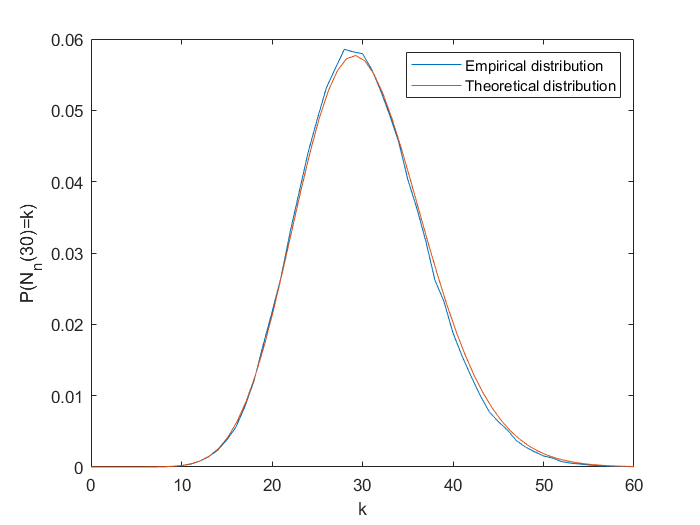}
\caption{Comparison between the empirical distributions of the number of visits of the observable $f(x)=(x-1/2)(x-1/4)$ in a ball centered at $0$ and the theoretical distribution described in the text for the map $3x \mod 1$.}
\label{vijj}
\end{figure}

We now give an example. We take the map $Tx=3x \mod1$, the observable $f(x)=(x-1/2)(x-1/4)$ and $f_0=0$. The two pre-images of $f_0$ are $1/2$ and $1/4$, of periods 1 and 2 respectively. From proposition 2, $\theta=7/9$. Then we have: $b_1=1/3$, $b_2=1/9$, $\mu_1=\mu_2=1/2$, $a_1=2t/9$ and $a_2=32t/81$. In figure \ref{vijj} we show the empirical distribution of the number of visits of $10^5$ different trajectories of length $10^6$ of the observable $f$ in the interval $I_r=(f_0-r,f_0+r)$, where $r=e^{-u}$, $u$ being the $0.995$-quantile of the distribution of $\phi$. We notice very good agreement with the theory.\\

\begin{figure}[t]
\centering
\includegraphics[height=4in]{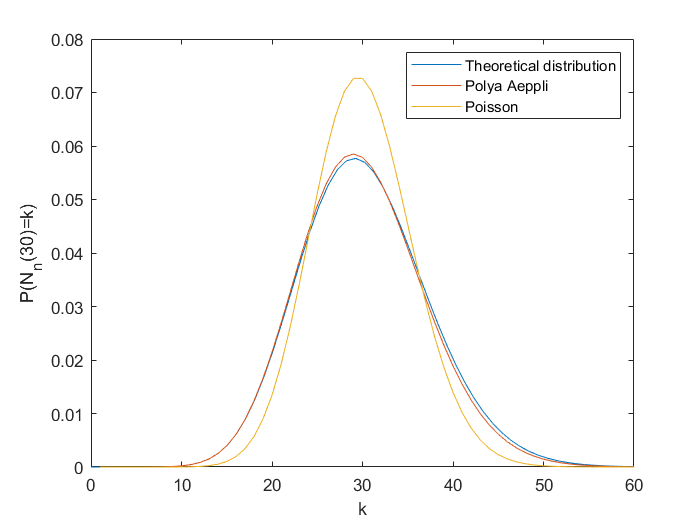}
\caption{Comparison between the distributions of the number of visits for the example in the text, a pure Poissonian distribution and a Polya Aeppli distribution with parameters given by $t=30$ and the extremal index $\theta=7/9$.}
\label{comp1}
\end{figure}

It is interesting to observe that if we take $\mu_1\neq \mu_2,$ but $b_1=b_2=b,$ which means we take the same periodicity for the two points $w_1, w_2,$ we recover the P\`olya-Aeppli distribution since $\lambda_l=b^{l-1}(1-b),$ with the extremal index $\theta=1-b.$ In fact,
even when $b_1\neq b_2$, numerical experiments suggest that the distribution stays close to a P\`olya-Aepply distribution. In figure \ref{comp1}, we show this effect by comparing the distribution associated with the example described in the text to a P\'olya-Aeppli distribution of parameters given by $\theta=7/9$ and $t=30$. The vicinity between the two distributions is striking and is found in a whole variety of examples. In \cite{ei}, we also observed this phenomenon in cases when the clustering structure is even more complex and for systems perturbed with discrete noise.

\begin{figure}[t]
\centering
\includegraphics[height=4in]{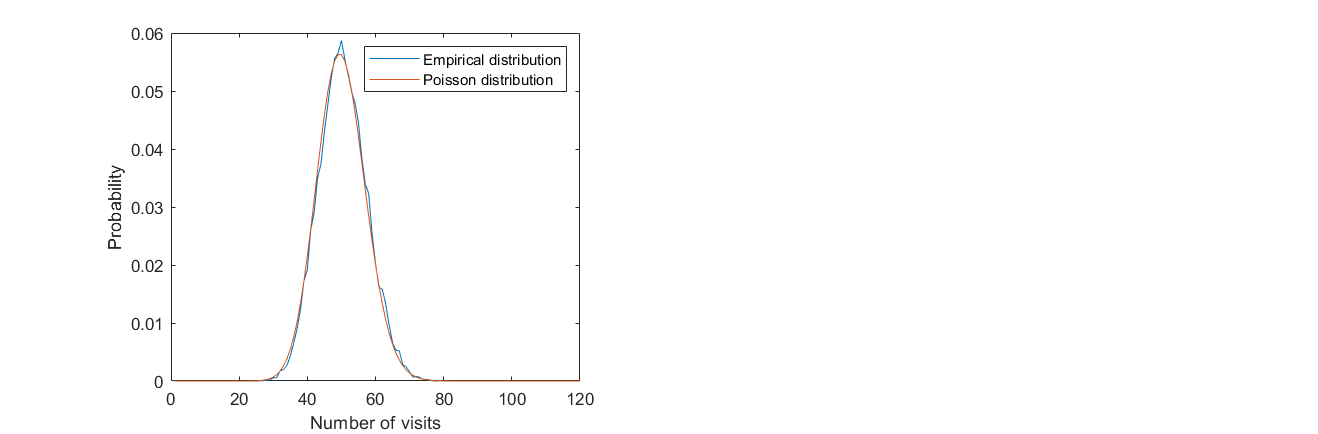}
\caption{Comparison between the empirical distributions of the number of visits of the observable ``mean value'' in a ball centered at $f(z)$ and a Poisson distribution for the baker map.}
\label{visbak1}
\end{figure}

 As we mentioned earlier, for a whole variety of observables $f$, no clustering is detected and the EI is $1$. We therefore expect to have a Poisson distribution for the statistics of visits. This is indeed what we observed for the baker map, with parameter $\alpha=1/3$, and the observable $f(x,y)=\frac{x+y}2$ (see figure \ref{visbak1}). We took a point $z$ at random in the attractor (actually we iterated a point in the basin several time to get it very close to the attractor),  and computed the empirical distribution of the number of visits of $10^5$ different trajectories of length $10^6$ for the observable $f$ in the interval $I_r=(f_0-r,f_0+r)$, where $r=e^{-u}$, $u$ being the $0.995$-quantile of the distribution.\\
 
\section{Randomly perturbed systems}
One could wonder what happens to the theory developed above when the dynamical system
is randomly perturbed; this has of course important physical applications when the system
    or its environment are affected by noise or when the available time series give only a partial
description of the evolution of the system variables.
As we anticipated in the Introduction, we will show that with suitable but very general choices
of the perturbations on the map or on the observable $f$ with values in $R^m$, the dimension of the
image measure will increase to $m$ if less than the dimension $m$ before perturbation, and drops
to $m$ otherwise.
\subsection{Perturbing the map}
  We defer to our paper \cite{ei} for an exhaustive presentation of different random perturbations in connection with  recurrence properties. For the purposes of this paper, we will consider random transformations, where  the iteration of the single map $T$ is replaced by the concatenation $T^n_{\underline{\omega}}:=T_{\omega_n}\circ \cdots \circ T_{\omega_1},$  where the $\omega_k\in \underline{\omega}:=(\omega_1, \cdots, \omega_k,\cdots) $ are i.i.d. random variables with (common) distribution $\mathbb{G}.$ Sometimes it is possible to show the existence of the so-called {\em stationary measure} $\rho_s$, verifying for any real bounded function $q$: $\int q d\rho_s=\int q\circ T_{\omega}d\rho_sd\mathbb{G},$ see \cite{book}
Chap. 7, for a general introduction to the matter. The product  $\mathbb{P}:=\rho_s\times \mathbb{G}^{\mathbb{N}}$ will give a stationary measure for the random process $q(T^n_{\underline{\omega}}(x), \sigma^n(\underline{w})),$ where $\sigma$ denotes the shift. The measure $\mathbb{P}$ will allow us to consider the limit theorems for such random processes in the so-called {\em annealed setting}; it will also weight the sets $B_{n,z}, C_{n,z}$ entering in the definition of the quantities $q_{n,k}$ expressing the extremal index. We defer to our papers \cite{ei} and \cite{Ay} for the analytic derivation of the extremal index in the annealed setting. We showed there in several examples, that whenever the distribution $\mathbb{G}$ has a density, the EI becomes equal to $1$, while it could be less than one for discrete distributions. The same happens in the present situation as the following two relatively simple situations indicate.
\begin{itemize}
\item {\em Continuous noise}. We consider a map $T$ verifying the assumptions in Proposition 2 and in particular we define it on the circle; we will say later how to generalize the result to the interval. We perturb $T$ with additive noise, namely we put $T_{\omega}(x)=T(x)+\omega$- mod $1$ and we choose $\omega$ with some smooth distribution $\mathbb{G}$ with density $q$ bounded from below. It is therefore possible to prove the existence of a stationary measure $\rho_s=h_s d\text{Leb}$ absolutely continuous with respect to Lebesgue with density $h_s.$ The computation of the extremal index follows now exactly the proof of Proposition 5.3 in \cite{Ay} with one difference: the connected ball $U_m$ there is now replaced by the set $f^{-1}B_{n,z}$ which is, in general, the disjoint union of a finite number of preimages. These sets are "centered" at the pre-images $\{z_l\}_{l\ge 1}$ of the target point $f(z).$ The key idea in \cite{Ay} was to show that for the majority of the realizations, with respect to $\mathbb{G}^{\mathbb{N}},$ the numerator in the quantities $q_{k,n}$ was zero. The rest was of higher order with respect to the denominator and vanished in the limit of large $n.$ The control in the numerator was based on the possibility to achieve, for a big portion of realizations $\underline{\omega},$ that $|T_{\underline{\omega}}^j(z)-z|>2 (\max{|T'|})^j |U_m|$, where $|U_m|$ denotes the diameter of $|U_m|$ and the latter is centered at $z.$ It easy to see that the same lower bound persists when the random orbit  $T_{\underline{\omega}}^j$ is computed starting from, say, $z_{l_1}$ and the right-hand side of the bound is replaced by the set $f^{-1}B_{n,z_{l_2}}$ around another point of the sequence $\{z_l\}_{l\ge 1}$. This is possible since the diameters of the $f^{-1}B_{n,z_l}$ are comparable, since $f$ is piece-wise $C^1.$ We left the details to the reader. At the end we get that all the $q_{l}=0,$ and therefore the extremal index $\theta=1.$ As we said above the proof extends easily to the additive perturbation of a piece-wise expanding map with finitely many branches verifying the other assumptions of Proposition 2.
    \item {\em Discrete noise} The purpose here is to give paradigmatic examples of the applicability of our theory with observable, leaving specific cases to other occasions.  For the discrete noise we could adapt to our first example  described at the beginning of section 3 with an invertible $f,$  the example studied in section 4.1.2 in our paper \cite{ei}. We considered there two maps on the circle $T_0=2x$- mod $1$ and $T_1=2x+b$- mod $1$, $0<b<1.$ If we now take the observable $f$ which is zero in $0$, $f(0)=0,$ we can repeat the argument in \cite{ei} with the set $B(0, e^{-u_n})$ there replaced by our $f^{-1}B(0, e^{-u_n}).$ The conclusion is that $q_0>0,$ and that $\theta<1.$
\end{itemize}
  We  argued in section 3 that in presence of observables the EI is difficult to compute; we believe that if in addition the system is randomly perturbed the EI is even more complicated and in general it should be $1$ or close to it. \\

  The computation of $\D$ in presence of noise is also interesting. We  first point out that our Proposition 1 easily generalizes to the annealed situation as we proved in \cite{ei} for discrete distributions  and in \cite{Ay} for distributions with density. Moreover,  we suppose that the target set is fixed and the  parameter $\tau$ defining the boundary level in Eq. (\ref{SCA3}) is independent of the noise, so that what we estimate {\em via} the convergence to the Gumbel law is the stationary measure of sets of type
$
\rho_s(B(f(z), \epsilon).
$
It is therefore interesting to evaluate that stationary measures; there are several ways to determine the existence of a stationary measures in connection with a given random perturbation, see for instance \cite{AA, viana, BHV}. Usually one needs a precise description of the stationary measures in order to establish stochastic stability, namely to recover the statistical properties of the unperturbed system when the noise is sent to zero. We are not interested in it; instead we are interested in getting an experimental way to construct a stationary measure and check its general properties. A useful result by Alves and Araujo \cite{AA} will provide us with   what we need. The idea is to look for a composition of maps close to a given one $T$ and  assume that the noise will verify two {\em nondegeneracy conditions}, namely:
\begin{itemize}
\item  (N1) The measure $\mathbb{G}^\mathbb{N}$ will be supported on a small set $S_{\epsilon}:=(\text{supp} \mathbb{G})^\mathbb{N}$ such that there is $\chi=\chi(\epsilon),$ for which each random orbit $T^n_{\underline{\omega}}(x)$ contains the ball of radius $\chi$ around $T^n(x)$ for all $x\in  X$ and $n$ sufficiently large.  As is written in \cite{AA}, this "condition means that perturbed iterates cover a full neighborhood of the unperturbed ones after a threshold for all sufficiently small noise level."
    \item (N2) We require that the measure $\int {\bf 1}_A(T^n_{\underline{\omega}}(x))d\mathbb{G}^{\mathbb{N}}
        (\underline{\omega})$ for any Borel set $A$  be absolutely continuous with respect to the Lebesgue measure $\text{Leb}$ on $X$, for all $x\in X$ and $n$ sufficiently large. This means that "sets of perturbation vectors of positive $ \mathbb{G}^\mathbb{N}$ measure must send any point $x\in X$ onto subsets of $X$ with positive Lebesgue measure after a finite number of iterates", \cite{AA}.
\end{itemize}
We now fix $x\in X$ and consider the measure, for any Borel set $A\subset X$:
\begin{equation}\label{RREE}
\rho_n(A):= \frac1n \sum_{j=0}^{n-1}\int {\bf 1}_A(T^j_{\underline{\omega}}(x))d\mathbb{G}^{\mathbb{N}}
        (\underline{\omega}).
\end{equation}
It has been proved in \cite{AA}, Lemma 3.5, that every weak* accumulation point of the sequence $\rho_n$ is stationary and absolutely continuous with respect to the Lebesgue measure whenever (N2) holds.
Notice that the Cesaro mean in Eq. (\ref{RREE}) is exactly the numerical procedure  to get the measure of a set by averaging over different realizations $\underline{\omega}$, so that we expect that with noise verifying the assumptions (N1) and (N2), the stationary measure $\rho_s$ is absolutely continuous with respect to Lebesgue \footnote{We notice that by general results on random perturbations, see for instance \cite{AA, BHV}, if the map $T$ preserves a unique absolutely continuous invariant measure, the absolutely continuous stationary measure is also unique.}.
This has an interesting consequence for the computation of $d_{\rho_s}^f$, since in presence of smooth observable $f$ and absolutely continuous measure $\rho_s$, Theorem 9 in \cite{JR} states that the dimension of the observable exists $\rho_s$-almost everywhere, is integer and is equal to the rank of $Df$.
We therefore expect that for such noises, the non-integer dimensions computed in the preceding examples for non-prevalent observable, become integer. For prevalent observable with large  dimensionality, $\D=D_1<m$, where $m$ is the dimension of the range of $f$, we expect that $\D$ drops to $\min(m,n)$, ($n$ being the dimension of the ambient space of the original system) in presence of a smooth stationary measure.\\

We tested this result by considering the dynamics on the product of two Cantor sets, with the non prevalent observable $f(x,y)=x-y$, and $f(z)=f(0,0)=0.$ The original dynamics given by an iterative function system is perturbed by an additive noise drawn with a uniform distribution in $B(0,\eta)$, for a small $\eta>0.$ To avoid that the dynamics leaves the unit square, we apply the mod-$1$ folding  after having applied the additive  perturbation. We observe in figure \ref{nois} that $\D$, which is about 0.63 when $\eta=0$ goes to 1 as $\eta$ increases. To compute $\D$, we simulated trajectories of $10^7$ points and considered blocks of size $10^3$.

Of course, if the noise does not verify assumptions (N1) and (N2), we do not know anymore if any weak limit of Eq. (\ref{RREE}) is absolutely continuous. This is in particular true if the unperturbed map $T$ will not preserve an absolutely continuous invariant measure. Otherwise and for uniformly expanding maps, it is always possible to get stationary measure which are absolutely continuous and that independently of the nature of the noise, \cite{BHV}.\\

If the stationary measure exists, the Hunt-Kaloshin Theorem still applies for the perturbed system, whatever the perturbation is. Indeed, this Theorem concerns measures and not the underlying dynamics.\\

When the perturbation is discrete, and the original measure has a fractal structure, the shape of the stationary measure is not yet completely understood. We therefore choose to study it numerically. We considered  the successive iterations of a baker map with $\lambda_a=\lambda_b=0.4$ and the parameter $\alpha$  equal to  $1/4$ and $1/3,$ each one with probability 1/2. For the observable $f(x,y)=(x+y,y^2),$ we found values for $\D$  around $1.70$ (we averaged the results over $20$ points of the attractor), which we interpreted as the local dimensions of the stationary measure. In fact, when we compute directly the local dimensions of this system, we also find a value of around 1.70. If we now we take a scalar observable $f(x,y)=x^2-y$, we find as expected values close to 1 for $\D$.

\subsection{Perturbing the observable}
We now suppose that the map $T$ does not change, but the observable does. In particular we assume that it changes in an i.i.d. way at each iteration. This could  have physical importance since it models  the influence of a random environment on the deterministic dynamics, or the uncertainty associated with the measurement process. By using the notations of section 1, we now consider the maximum of the random  variables, for $k=0,\cdots, n-1:$
$$
\phi(T^k(x),\omega_k)=-\log(\text{dist}(f_{\omega_k}(T^k(x))-f(z))),
$$
where the $\omega_k$ are i.i.d. random variable with common distribution $\mathbb{G}$ and $f(z)$  is the value of a fixed observable at the point $z.$ The probability will now be  $\mu\times \mathbb{G}^{\mathbb{N}}$  and we indicate it with $\mathbb{P}.$ We write again $\underline{\omega}$ for the vector with components $\{\omega_k\}_{k=0\cdots, \infty}.$ The maximum will therefore be a function of $x$ and $\underline{\omega}$, $M_n:=M_n(x, \underline{\omega}).$  By setting ourselves in the framework of the uniformly expanding maps considered in section 1,  we immediately have
$$
\mathbb{P}(M_n\le u_n)=\int ({\bf 1}_{C_{n,z}}\circ f_{\omega_0})(x) \cdots ({\bf 1}_{C_{n,z}}\circ f_{\omega_{n-1}})(T^{n-1}x) h(x) dx d\mathbb{G}^{\mathbb{N}}(\underline{\omega}),
$$
where the $C_{n,z}$ have the same meaning as in section 1. Since the $\{\omega_k\}_{k=0\cdots, \infty}$ are independent and performing first the integration with respect to $\mathbb{G}^{\mathbb{N}}$, we have
$$
\mathbb{P}(M_n\le u_n)=\int U_n(x) U(Tx)\cdots U_n(T^{n-1}x)h(x)dx,
$$
where
$$
U_n(x):= \int ({\bf 1}_{C_{n,z}}\circ f_{\omega})(x)d\mathbb{G}(\omega).
$$
For instance, if we keep an initial $f$ with value in $\mathbb{R}$ and add to it a random term $\eta$ with uniform distribution in $[-a, a]$ we have
$$
U_n(x)= \frac{1}{2a}\int_{-a}^a {\bf 1}_{C_{n,z}}(f(x)+\eta) d\eta=\frac{1}{2a}\text{Leb}\{ [-a,a]\cap [C_{n,z}-f(x)]\}.
$$
Another choice is to add to an unperturbed observable $f$ two quantities $\eta_1, \eta_2$ taken with respective probabilities $p_1, p_2.$ In this case we have
$$
U_n(x)= {\bf 1}_{f^{-1}[C_{n,z}-\eta_1]}(x)p_1+{\bf 1}_{f^{-1}[C_{n,z}-\eta_2]}(x)p_2.
$$
Then
$$
\mathbb{P}(M_n\le u_n)=\int \tilde{P}_n^n(h)(x) dx,
$$
where $\tilde{P}_n(h):=P(U_nh).$\\

We are now in position to apply the spectral theory since we have just constructed a REPFO system: we leave the details to the reader in order to check the necessary requirements. What is important for us now, is to give an expression for the extremal index and for the boundary level $u_n$, which will reflect on the dimension of the image of the observable.  Let us begin with the extremal index. The quantities $q_{k,n}$ are now defined as \cite{GK, KL}:
$$
q_{k,n}=\frac{\int (P-\tilde{P}_n)\tilde{P}_n^k(P-\tilde{P}_n)(h)dx}{\int (P-\tilde{P}_n)(h)dx}.
$$
By posing
$$
V_n(x):= \int ({\bf 1}_{B_{n,z}}\circ f_{\omega})(x)d\mathbb{G}(\omega),
$$
we immediately have
$$
q_{k,n}=\frac{\int V_n(T^{k+1}(x))U_n(T^k(x))\cdots U_n(T(x))V_n(x)d\mu}{\int V_n d\mu},
$$
which allows us to construct the EI $\theta.$\\ As in the previous section, we now give the computation of the EI in two situations, with continuous and discrete noise.
\begin{itemize}
\item {\em Continuous noise}

We put ourselves in the setting of Proposition 2 plus other assumptions which we will add during the proof. Let us consider the additive  noise $f(x)+\eta$ described above with $\eta$ much smaller than $1$. The first and the last terms in the integral in the numerator of the $q_{k,n}$ are:
  \begin{equation}\label{ca}
    \frac{1}{2a}\int_{-a}^a {\bf 1}_{B_{n,z}}(f(x)+\eta) d\eta \ \text{and} \ \frac{1}{2a}\int_{-a}^a {\bf 1}_{B_{n,z}}(f(T^{k+1}x)+\eta) d\eta.
    \end{equation}
    In particular both quantities are bounded by 
    $$
    \frac{1}{2a}\text{Leb}[(B_{n,z}-f(T^{j}x))\cap [-a,a]]\le\frac{1}{2a} \text{Leb}(B_{n,z}), \ j=0, k+1,
    $$
    and therefore the numerator in $q_{k,n}$ is bounded from above by $\frac{1}{4a^2} (\text{Leb}(B_{n,z}))^2.$ \\
    We now rewrite the denominator as
    $$
    \int  \frac{1}{2a}\int_{-a}^a {\bf 1}_{B_{n,z}}(f(x)+\eta) d\eta d\mu= \frac{1}{2a}\int_{-a}^a\left[\mu(f^{-1}(B_{n,z}-\eta))\right]d\eta.
    $$
    We now suppose that the preimages of the set $B_{n,z}-\eta$ are at most $L$ for any $\eta$ and set $\max |f'|=M_D<\infty$; moreover we suppose that the density $h$ of $\mu$ is bounded from below by $h_m$. Then
    $$
    \mu(f^{-1}(B_{n,z}-\eta))\ge Lh_m M_D^{-1}\text{Leb}(B_{n,z}-\eta)
    $$
    which implies after integration with respect to $\eta$:
    $$
     \int  \frac{1}{2a}\int_{-a}^a {\bf 1}_{B_{n,z}}(f(x)+\eta) d\eta d\mu\ge \frac12Lh_m M_D^{-1}\text{Leb}(B_{n,z}).
    $$
    If we now divide the numerator with the denominator, we will find the ratio going to zero for $n\rightarrow \infty$, which shows that all the $q_{k,n}$ are zero and the EI is one.
\item{\em Discrete noise}
    
  We give this  example again   in the setting of Proposition 2. Take the discrete noise with distributions $\{(\eta_1, p_1), (\eta_2, p_2)\},$  and the ball $B_{n,z}$ around the point $f(z).$ Put $B_{n,z,1}=B_{n,z}-\eta_1, B_{n,z,2}=B_{n,z}-\eta_2$ and $\eta_1<0<\eta_2.$ If $z$ is a point where $f$ is monotone and we choose $n$ sufficiently large,  the sets  $B_{n,z,1}, B_{n,z,2}$ will be disjoint and the same for the four sets $f^{-1}_i B_{n,z,j}, i,j=1,2.$ Call $z_{i,j}$ the point such that $f(z_{i,j})\in B_{n,z,j}, f(z_{i,j})+\eta_j=f(z)$. Suppose now that the point $z_{1,1}$ is a fixed point for $T,$ but the remaining points $z_{i,j}$ are not periodic for $T$. Then the only term which could give a non zero contribution is $q_{0,n}$, which reads
$$
q_{0,n}=\frac{\int [{\bf 1}_{f^{-1}[B_{n,z,1}]}(x)p_1+{\bf 1}_{f^{-1}[B_{n,z,2}]}(x)p_2][{\bf 1}_{f^{-1}[B_{n,z,1}]}(Tx)p_1+{\bf 1}_{f^{-1}[B_{n,z,2}]}(Tx)p_2]h(x)dx}{\int[{\bf 1}_{f^{-1}[B_{n,z,1}]}(x)p_1+{\bf 1}_{f^{-1}[B_{n,z,2}]}(x)p_2]h(x)dx}.
$$
When $n$ goes to infinity only the term
$$
p_1^2 \mu(f^{-1}_1(B_{n,z,1})\cap T^{-1}f^{-1}_1(B_{n,z,1}))
$$
gives a non zero contribution. By using the same distortion arguments as in the proof of Proposition 2 we immediately get
$$
q_{0}=\frac{p_1}{|T'(z_{1,1})|}\left[ 1+\sum_{i,j=1,2; i\neq j}\frac{p_j}{p_1}\frac{|f'(z_{1,1})|h(z_{i,j})}{|f'(z_{i,j})|h(z_{1,1})}\right]^{-1}
$$
and the extremal index will be $\theta=1-q_{0}$.
\end{itemize}

We now discuss the choice of the boundary levels $u_n.$ First it is defined as
$$
\mathbb{P}\left((x, \omega); -\log|f_{\omega}(x)-f(z)|>u_n\right)\rightarrow \tau/n.
$$
By introducing the image measures
$$
\mu_{\omega}^*:=f_{\omega}^*\mu,
$$
we finally have
$$
\int \mu_{\omega}^* (B(f(z), e^{-u_n})d\mathbb{G}(\omega)\rightarrow \tau/n.
$$
It is  interesting to explore whether we have a scaling of type
$$
\mu_{\omega}^* (B(f(z), r)\approx r^{d^*},
$$
and finally
$$
\int \mu_{\omega}^* (B(f(z),r)d\mathbb{G}(\omega)\approx \int r^{d^*}d\mathbb{G}(\omega)
$$
for some exponent $d^*.$
Notice that contrarily to formula (\ref{LLL}) we are now transporting the measure $\mu$ with some $f_{\omega}$ and computing this  measure around the image of a point with a different $f$.

A  simple trick  allows us to restore the right framework and a quite general example will suggest some expected behavior. Take a countable family of prevalent observable
indexed by $f_j, j=1, \dots, \infty$ each with a weight $p_j$ such that $\sum_j p_j=1.$  This discrete measure is called $\mathbb{G}$. Fix one $f$ and suppose that the range of each $f_j$ contains $f(x)$; set $\mu_j^*=f_j^*\mu.$ Then
$$
\mu_{j}^* (B(f(x), r))=\mu\left(f_{j}^{-1}(B(f_{j}(f^{-1}_{j}(f(x))), r))\right).
$$
Call $x_{j}$ one of the pre-images of $f(x)$ by $f_{j}$, $x_{j}\in f_{j}^{-1}(f(x))$. Then
$$
\mu_{j}^* (B(f(x), r))=\mu_{j}^* (B(f_{j}(x_{j}), r)).
$$
Each $f_j$ is prevalent  so by Theorem 1 we know that the quantity $\underline{d}^{f_j}_{\mu}(z)$ is equal to the minimum between the dimension of the range of $f,$ which we take  equal to $m,$ and the lower point-wise dimension of $\mu$ at $z$, provided the latter is chosen  $\mu$-a.e. If we suppose that the point $x_j$ is typical for $f_j$ and also that $\mu$ is exact dimensional, we have that
$$
\mu_{j}^* (B(f(x), r)\approx r^{d^{f_j}_{\mu}(x_j)}, \text{with} \ d^{f_j}_{\mu}(x_j)=\min(m, D_1(\mu)).
$$
In conclusion
$$
\int \mu_{\omega}^* (B(f(z),r)d\mathbb{G}(\omega)\approx
r^{(\min(m, D_1(\mu)))}.
$$
Therefore for scalar functions $(m=1)$ and attractors of high dimensionality, we expect to get a dimension equal to $1$ when the observable is perturbed. Instead if the attractor, or repeller, have dimension less than $m$, the dimension of the image measure  will jump to $m.$ 

\begin{figure}
     \begin{subfigure}[t]{0.5\textwidth}
        \centering
\includegraphics[height=2.5in]{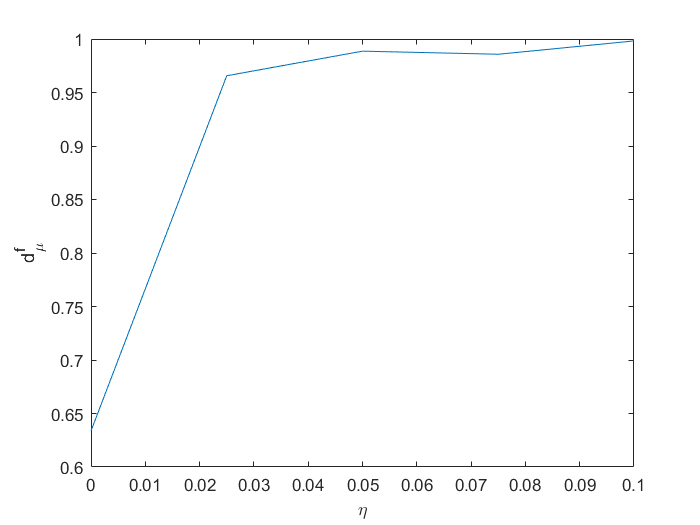}
\label{visbak}
    \end{subfigure}%
    ~
    \begin{subfigure}[t]{0.5\textwidth}
        \centering
     \includegraphics[height=2.5in]{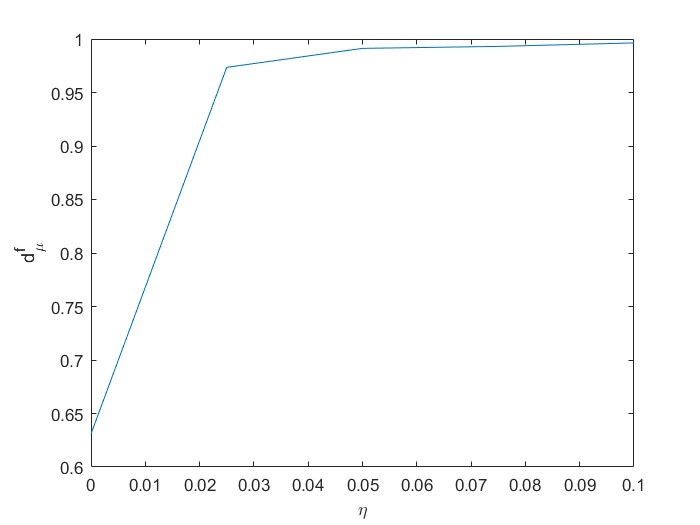}
     \label{d}
    \end{subfigure}
    \caption{Influence of uniform noise of different intensities applied to the dynamics (left) and to the observable (right), for the motion on the product of Cantor sets and the non-prevalent observable $f(x,y)=x-y$, with $f_0=f(0,0)=0.$}
    \label{nois}
\end{figure}

We studied the effect of uniform additive noise of different intensities to the observable $f(x,y)=x-y$ for the dynamics on the product of two  Cantor sets described earlier. At each iteration, we computed $f_k=f(T^k(x,y))+\varepsilon_k$, where $\varepsilon_k$ are i.i.d. random variables drawn with a uniform distribution in $[-\eta, \eta].$ Results are shown in figure \ref{nois}. Similarly to the case where the dynamics is perturbed, we observe a convergence of $\D$ to 1 as the intensity of noise $\eta$ increases. We stress that this monotonic convergence to $1$ depicted in the figure is a numerical artefact, since the image dimension becomes immediately $1$ as soon as the noise is switched on. To compute $\D$, we simulated trajectories of $10^7$ points and considered blocks of  $10^3$ points.

\section{Open systems}
In the paper \cite{targets} we considered the extreme value distribution for open systems, namely for systems with holes, where the orbits enter and disappear forever.  That was motivated by  the statistical description of phenomena where a perishable dynamics
 is approaching a fixed target state,  but at the same time it deviates to another  location where it is captured or vanishes.  It is useful to extend that theory in presence of observables. A close look at the proofs in the aforementioned paper, shows that such proofs can be easily translated to our present situations. One of the major results in \cite{targets} was to relate the extremal index to the escape rate (from the hole). That was achieved when the target set was chosen around periodic points. In presence of observables, periodicity is much more cumbersome, as we described in Proposition 2; it would be therefore interesting to have a version of such a proposition  in the presence of holes. Before doing that we  recall the main result in \cite{targets}.
\begin{proposition} \cite{targets}\label{FRT}
Let $T$ be a uniformly expanding map of the interval $I$ preserving a mixing
measure. Let us fix a small absorbing region, a hole $H\subset I$; then there
is an absolutely continuous conditionally invariant measure $\nu,$
supported on $X_0=I\setminus H$ with density~$h_0.$ Write
$\alpha=\nu(T^{-1}X_0)$.  If the hole is small enough there is a
probability measure $\mu_0$ supported on the surviving set $X_{\infty}$ such that
the measure $\Lambda=h_0\mu_0$ is $T$-invariant; we  assume that $h_0$ is bounded
away from zero. Having fixed the positive number $\tau$, we take the sequence $u_n$ satisfying $n\Lambda(B(z,\exp(-u_n) ))=\tau,$ where~$z\in X_{\infty}$.  Then, we take the sequence of conditional probability measures $\mathbb{P}_n(A)=\frac{\nu(A\cap X_{n-1})}{\nu(X_{n-1})},$ for $A \subset I$ measurable, and define the random variable $M_n (x):=\max\{\phi(x), \dots,
\phi(T^{n-1}x)\},$ where $\phi(x)=-\log|x-z|.$ Moreover we  suppose that all the iterates $T^n, n\ge 1$ are continuous at $z$ and also that $h_0$ is continuous
at $z$ when the latter is a  periodic point. Then we have:
\begin{itemize}
\item If $z$ is not a periodic point:
\[
\mathbb{P}_n(M_n\le u_n)\rightarrow e^{-\tau}.
\]
\item If $z$ is a periodic point of minimal period $p$, then
\[
\mathbb{P}_n(M_n\le u_n)\rightarrow e^{-\tau\theta},
\]
where the extremal index $\theta$ is given by:
\[
\theta=1-\frac{1}{\alpha^{p}|(T^p)'|(z)}.
\]
\end{itemize}
\end{proposition}
We remind that a probability measure $\nu$ which is absolutely continuous with respect to Lebesgue is called {\em a conditionally invariant  probability measure} if it satisfies for any Borel set $A\subset I$  and for all $n>0$ that
\begin{equation} \label{eq:cipm}
\nu(T^{-n}A\cap X_n)= \nu(A)\ \nu(X_n).
\end{equation}
 The \emph{surviving set} is defined as $\smash{X_{\infty} = \bigcap_{n=1}^{\infty} X_n}$, where $\smash{X_n=\bigcap_{i=0}^n T^{-i}X_0}$ is the set of points that have not yet fallen into the hole at time~$n$. Finally   the {\em escape rate} $\eta$  for our
open system is usually defined as $ \eta=-\log \alpha.$

Let us return to the proof of  Proposition 2 trying to adapt it. The class of maps are the same as those in  Proposition \ref{FRT}. The main change will concern the invariant measure which is now the singular measure $\Lambda$ on the surviving (fractal) invariant set $
X_{\infty}.$ Such an invariant measure is absolutely continuous with respect to the conformal measure called $\mu_0$ in Proposition \ref{FRT}. This conformal measure plays the role of the Lebesgue measure in the proof of Proposition 2; in the latter we performed a change of variable which produced the terms $|(T^p)'|$, where $p$ was related to the periodicity of the point where we computed the derivative. The conformal measure will give a multiplicative factor $\alpha^p.$ Moreover the density $h$ in Proposition 2 will be now replaced with the density $h_0$ with respect to the conformal measure $\mu_0.$ In conclusion the term $q_k$ in (\ref{gen}) will be now replaced by the following one, which we call $q_k^{(o)}$ since it refers to open systems
\begin{equation}\label{gen2}
q_k^{(o)}=\underset{w\in A_k}\sum \frac{1}{\alpha^{k+1}|T^{(k+1)}(w)'|}\frac{1}{1+ \frac{|f'(w)|}{h_0(w)}\underset{y\in B_k}\sum \frac{h_0(y)}{|f'(y)|}}.
\end{equation}
If we want to perform numerical computations, we should know the value of $\alpha.$ We already said that $\alpha$ is related to the size of the hole, in particular one can show that $\alpha$ is the largest eigenvalue of the perturbed transfer operator $\hat{P}g:=P({\bf 1}_{H^c} g)$, compare with the perturbed operator $\tilde{P}_n$ of section 2. Therefore for small hole one could apply again the spectral technique of \cite{KL} and get $\alpha$ as an asymptotic perturbation of $1$, the largest eigenvalue of $P.$ It is not therefore surprising that such an expansion will be related to the location of the hole. In particular if the latter is around a point $z$ which is not periodic, $\alpha$ will be equal to $1$, instead it will be equal to $1-\frac{1}{|(T^p)'|(z)}$ if $z$ is a periodic point of minimal period $p.$ We point out again that those values hold in the limit of vanishing holes, so that one would get something slightly different for hole with finite size. An interesting case of a large  hole is given in the next section.

\subsection{EVT on fractals I} 
In this section and in its companion  8.2, we address the following question. Suppose we have a fractal invariant set which is a repeller and whose Lebesgue measure is zero. How could we get a good extreme value theory by using the Lebesgue measure as the underlying probability? In fact almost all the orbits leaving on sets of positive Lebesgue measure tend to escape from the repeller.  On the other hand Lebesgue measure is the most accessible measure and repellers are widespread objects, for instance they constitute the basin boundaries between two, or more, basin of attraction, see \cite{BL} for applications to climate.  The simplest non-trivial repeller is probably the ternary Cantor set, $\mathcal{C}$; in the above terminology,  it is the   surviving set of the map $T(x)=3x$-mod$1$ having taken the hole as the open interval $(1/3, 2/3).$ We point out that other repellers  could be generated as the surviving  sets in open systems, so that the next considerations could be useful to understand larger class of fractal invariant sets.
The first study dealing with the ternary Cantor set in connection with EVT was mostly numerical and it was given in \cite{mant}: the authors conjectured the existence of a limiting extreme value law with an EI equal to $1$. A rigorous proof appeared recently in the paper \cite{rf}; in particular, the authors introduced the observable
$$
\phi(x)=\bigg\{
\begin{array}{rl}
n, & n\in C_n\\
\infty, & \text{otherwise}\\
\end{array}
$$
where $C_n$ is the disjoint union of the $2^n$ sets in the construction of the Cantor set $\mathcal{C}$.\footnote{The Cantor set is given by $\mathcal{C}=\cap_{n\ge 1}C_n,$ where the $C_n$ denotes the disjoint union of the $2^n$ (cylinder) sets obtained by removing the middle third part of each connected component of of $C_{n-1}.$}\\
Notice that the function $\phi$ will have his maximum (infinity) on the Cantor set, otherwise it says how fare we are from it: it is called the Cantor ladder function in \cite{mant, rf}. The probability was chosen as the Lebesgue measure $\text{Leb}$ on the unit interval. Given $\tau>0$ and by introducing the sequence of thresholds 
$$
w_n:= \left\lfloor{\tau\left(\frac{3}{2}\right)^n}\right\rfloor,
$$
it was proved in \cite{rf} that
$$
\lim_{n\rightarrow \infty} \text{Leb}\left(M_{w_n}\le n\right)=e^{-\tau(1-\frac23)},
$$
where $M_n$ is the process as defined in (\ref{MAX}). In this setting, the EI is therefore equal to $1/3.$ This result is interesting since the limit distribution is obtained with the Lebesgue measure, which  allows us to look at the whole Cantor set as a rare event. 

We now instead provide  a local inspection to the Cantor set by giving the statistics of the hitting time around any point on the repeller. This statistics will be given by a measure which is absolutely continuous with respect to Lebesgue. All this will follow automatically from our Theorem \ref{FRT} if it would hold  for such a big hole like $H=(1/3, 2/3).$
 Actually, in that  theorem  we required the hole to be small to be able to construct the conformal measure $\mu_0$ and its density $h_0$ with a perturbative argument.  In our case, we can do it directly since the map is easy enough. If we set $P$ the Perron-Frobenius operator associated to $T$ and we define the perturbed operator $P_0$ as $P_0(g)=P(g {\bf 1}_{H^c}),$ where $g$ is a function of bounded variation, we check easily that, having set $\alpha=\frac23$:
\begin{eqnarray}
P_0 h_0=\alpha h_0\\
P_0^* \mu_0=\alpha \mu_0,
\end{eqnarray}
  where:   $P_0^*$ is the dual of $P_0;$ $h_0=3/2$ on the unit interval and $\mu_0$ is the balanced measure described below.  The absolutely continuous conditionally invariant measure $\nu$ will be a measure with density $3/2$ on the closed intervals $[0, 1/3], [2/3, 1],$ and $0$ on the hole. Using it to construct the absolutely continuous probability $\mathbb{P}_n$ given in Theorem  \ref{FRT}, which is numerically accessible, we could place target sets around any point $z\in \mathcal{C},$  as balls of radius $e^{-u_n},$ where the thresholds $u_n$ can be chosen as
  $\frac32n\mu_0(B(z, e^{-u_n}))=\tau.$ Therefore we get convergence to  Gumbel's law for our process $M_n$ (\ref{MAX}) with:\\
  - the EI is equal to $1$ if the point $z$ is not periodic,  thus partially supporting the conclusions of \cite{mant}.\\
  - if we choose $z$ as a periodic point (they are dense in $\mathcal{C}$),  of minimal period $p,$  we get for the EI $\theta:$
  $$
  \theta= 1- \frac{1}{2^p}.
  $$
The global \cite{rf} and  our local  approaches to the EVT distribution just described,  considered the Cantor set as the non-wandering set of a dynamical system defined on the unit interval. One could consider the dynamics defined directly on the Cantor sets, which means to study the system $(\mathcal{C}, T_{|\mathcal{C}}).$ We took this point of view in \cite{D2}, where the transfer operator was defined directly  on $\mathcal{C}$ with the potential $|T'|^{-d_H},$ where $d_H$ was the Hausdorff dimension of the Cantor set ($d_H=\frac{\log 2}{\log 3}$ for the ternary Cantor set). It turns out that the invariant (Gibbs) measure for that potential is exactly the measure $\mu_0$ introduced above \footnote{Notice that $\mu_0$ is also invariant in the framework of Proposition \ref{FRT} since it differs from the measure $\Lambda$ by the constant $f_0=3/2.$ In this respect $\Lambda$ is not a  probability measure.}. But there was another reason for having chosen such a potential; in fact  the conformality of this measure  implies that for any measurable set $A$ where $T$ is one-to-one, we have $\mu_0(TA)=2 \mu_0(A).$ Therefore  that measure gives  $2^{-n}$ masses to the $2^n$ intervals $C_n$ of length $3^{-n}$ at the $n$-th generation in the construction of $\mathcal{C}.$  One could also show that this measure  is the weak-limit of the sequence of point masses measures  constructed with the pre-images of {\em each} point in the interval both weighted by $1/2.$ This is a sort of ergodic theorem for repellers, which makes $\mu_0$
  accessible  for numerical purposes: it is often called a {\em balanced} measure. {\em Using $\mu_0$ as the probability for the EVT distribution directly on the Cantor set, we find  Gumbel laws with the same behavior for the EI described above}. We will use again this balanced measure in section 8.2.

\section{Hitting time statistics in the neighborhood of sets}
 \subsection{Smooth sets}
 Our approach allows us to compute the hitting time statistics (and the statistics of the number of visits) in shrinking neighborhoods of a $C^1$ surface $\Gamma \subset \mathbb{R}^k$. At this regard, it is enough to consider an observable $f\in C^1(\mathbb{R}^k,\mathbb{R})$ such that $f(x)=\text{dist}(x,\Gamma)$. In this case, we have for all $n$ the identity
$$\{x\in X, |f(x)-0|<e^{-u_n}\}=\{x\in X, \text{dist}(x,\Gamma)<e^{-u_n}\}.$$

The hitting time statistics in the target sets $\Gamma_n =\{x\in X, |f(x)-0|<e^{-u_n}\}$ can be deduced from our theory and is given by the distribution of $M_n$, which converges to the Gumbel law. We then automatically obtain the hitting time statistics in the sets $\{x\in X, \text{dist}(x,\Gamma)<e^{-u_n}\}$. The parameters of this limit law are often computable explicitly (see section \ref{dddd} for the computation of $\D$).\\

We now give two examples based on the baker map and on the product of two Cantor sets.

\begin{itemize}
    \item Let us take $\Gamma$ a straight line of equation $ax+by+c=0$. The distance from a point $z=(x,y)\in X$ to $\Gamma$ is given by $$\text{dist}(z,\Gamma)=\frac{|ax+by+c|}{\sqrt{a^2+b^2}}.$$
    
    Let us take the $C^1$ observable
     $$f(x,y)=\frac{ax+by+c}{\sqrt{a^2+b^2}},$$
    
   so that with this choice of observable, we have for all $n$ the identity
    \begin{equation}\label{GRT}
   \{z\in X, \text{dist}(z,\Gamma)<e^{-u_n}\}=\{z\in X, |f(z)-0|<e^{-u_n}\}.
   \end{equation}
   
   The hitting times statistics in the set in the right hand side is given for large $n$ by the Gumbel law with scale parameter $1/\D$. For the baker's map $\D$ is $1$ if $b\neq 0$ and less than $1$ if $b=0$ (see section \ref{dddd}). The extremal index was computed numerically by two of us in \cite{D2} in the neighborhood of the diagonal and we found a value strictly less than $1.$  \\ For the product of the two Cantor sets, we found in section \ref{dddd} that for straight lines parallel to the coordinate axis and in the neighborhood of the diagonal, $\D$ was strictly less than one. In \cite{D2} we proved analytically that the EI computed around the diagonal was equal to $1/2.$ 
   
    \item Let us take now  $\Gamma$ as the circle in $\mathbb{R}^2$  of center $(a,b)$ and of radius $R$, of equation $$(x-a)^2+(y-b)^2-R^2=0.$$
    
    The distance from a point $z=(x,y)$ to $\Gamma$ is given by 
    $$\text{dist}(z,\Gamma)=|\sqrt{(x-a)^2+(y-b)^2}-R|.$$
    
    Let us take the $C^1$ observable 
     $$f(x,y)=\sqrt{(x-a)^2+(y-b)^2}-R,$$
    
     so that with this choice, we have again the equivalence (\ref{GRT}).
   
   
   As before, the hitting times statistics in the neighborhood of the circle is given for large $n$ by the Gumbel law with scale parameter $1/\D$. For the baker's map $\D$   is $D_1$ when $R=0$ as we already showed, and $1$  whenever $R>0,$ as it easy to see by adapting the argument given for the double Cantor set in the neighborhood of the diagonal. \\ For the product of two Cantor sets we have again $\D<1$ for $R=0$ and also $\D<1$ for $R>0,$  proving that the observable is not prevalent. In both cases the EI follows the usual dichotomy for $R=0$. We do not dispose of rigorous results in the other case $R>0.$
    
\end{itemize}

One can generalize this approach to higher dimensional $REPFO$ systems and generic $C^1$ hypersurfaces.

\subsection{EVT of fractals II} 
We could now wonder what happens if we consider the distance with respect to a fractal set, for instance the ternary Cantor set introduced in section \ref{dddd}. An easy way to do it and which uses the ideas of this section, is to consider the ternary Cantor set placed along the $y$-axis in the Cartesian product studied in section 4.2. If we consider the function $f(x,y)=x,$ we are led to study the EVT distribution for the observable (\ref{OB}) with $\phi(x,y)=-\log(|x|).$ As the underlying probability we take the balanced measure $\mu^{(2)}$ introduced in section 4.2 and described in section 7.1.  We are now interested in computing the EI. Since the Cantor set on the $y$-axis is invariant for the product map $T\times T$ on the unit square (the map $T$ was defined in section 4.2), we can adapt the proof "along the diagonal" given by us in \cite{D2} section II B or in \cite{D1, HV2},  and find easily that only the term $q_0$ in the expansion of the extremal index will not vanish. Then we use the conformality of the factor measure  $\mu$ along the $x$-axis and we  get   $q_0=1/2,$ giving also an extremal index equal to $0.5.$ It is interesting to compare this result with that in \cite{rf} described in section 7.1 and with a global approach to the Cantor set as a rare event: the EV found there was $1/3.$\\

Up to now, EVT has been  applied to compute hitting time statistics in the neighborhood of some sets of points \cite{FFT,uncount}, some Cantor sets \cite{rf,giorgiofractal}, or the diagonal in product spaces \cite{D2,dq}. We now provide generalizations  to arbitrary $C^1$ surfaces. We point out that a few results have been obtained in that direction in \cite{nicol}, where for the Arnol'd cat map and a $C^1$ curve $\Gamma \subset [0,1]^2$, the asymptotic behavior of the shortest distance of the system to $\Gamma$ up to a time $n$ was derived for Lebesgue almost every starting point.

\section{Large deviations}
 We pointed out in the Introduction and experienced in the preceding sections,  that one the most useful, and used, consequences of the EVT applied to dynamical systems, is the possibility to compute numerically the point-wise (also named local),  dimensions of the invariant sets. It turns out that in several time series given by natural phenomena or experimental signals, these local dimensions deviate significantly from each other, while in the ergodic setting they should coincide almost everywhere. Instead of seeing in this behavior only a numerical effect, we attributed it to the presence of large deviations in the convergence to the local dimension. The latter manifest  themselves on small, but not negligible, scales, a regime which we called {\em penultimate} \cite{dq}.  The presence of large deviations  for the point-wise dimensions has been rigorously proved for conformal repellers\footnote{These are the invariant sets of uniformly expanding $C^{1+\alpha}$ maps, defined
on smooth manifolds and whose derivative is a scalar times an isometry. The
repeller arises as the attractor of pre-images of the map, see \cite{21} for an
exhaustive description. Dynamically generated Cantor sets on the line, Iterated
Function Systems with the open set condition, disconnected hyperbolic Julia sets,
are all examples \cite{22} of conformal repellers. It is worth mentioning that such
repellers can be coded by a subshift of a finite type and they support invariant
measures which are Gibbs equilibrium states. This makes them particularly
suited for the application of the thermodynamic formalism.} in the paper \cite{gdlocdim}. Suppose we have an exact dimensional measure $\mu$, call $D_1$ the $\mu$-almost sure limit, and suppose that the following limit exists
 \begin{equation}\label{gendi}
 D_q=\lim_{r\to0}\frac{\log \int \mu(B(x,r))^{q-1}d f_*\mu(x)}{(q-1)\log r}.
 \end{equation}
 for all $q\in \mathbb{R}$ and moreover the function $\tau(q)=D_q(q-1)$ is $C^1$ over $\mathbb{R}$ and strictly convex.\footnote{For $q=1,$ the value for $D_q$ is obtained by l’Hopital rule.} This is what happens  for conformal repellers, where the limit (\ref{gendi}) exists for real  $q$ \cite{PW}. Then we are in the setting of the large deviation result by Gardner-Ellis, see for instance  \cite{ZD},  which allows us to state for all interval $I$:
  \begin{equation}\label{lhi}
\lim_{r\to0}\frac{1}{{\log r}} \log \mu\left(\left\{z \in X \mbox{ s.t. } \frac{\log\mu(B(z,r))}{\log r}\in I\right\}\right)=\inf_{s\in I}Q(s).
\end{equation}
 The {\em rate function} $Q(s)$ is determined by the $D_q$:
\begin{equation}
Q(s)=\sup_{q\in \mathbb{R}}\{-qs+qD_{q+1}\}.
\end{equation}
\begin{remark}\label{rid}
We notice that when the limit (\ref{gendi})  exists in some interval of values of $q,$ then we have to restrict   the interval $I$ to a suitable  neighborhood $\tilde{I}$ of the information dimension $D_1$, and for $s\in \tilde{I}$ we can control only deviations larger than $D_1,$ namely we have
\begin{equation}\label{secld}
\lim_{r\to0}\frac{1}{{\log r}} \log \mu\left(\left\{z \in X \mbox{ s.t. } \frac{\log\mu(B(z,r))}{\log r}>D_1+s\right\}\right)=\inf_{s\in \tilde{I}}Q(D_1+s).
\end{equation}
see \cite{HHHH}, Lemma XIII.2, for the details.  
\end{remark}

 It is interesting to ask whether large deviations are present when an observable is applied to the measure. Let us start by defining the generalized dimension of order $q$ of $f_*\mu$ (if it exists) as:
 \begin{equation}\label{mugedi}
 D_q^f=\lim_{r\to0}\frac{\log \int f_*\mu(B(x,r))^{q-1}d f_*\mu(x)}{(q-1)\log r}.
 \end{equation}
 Suppose now that the image measure $f_*\mu$ is exact dimensional; 
 if the function $D_q^f$ exists and is differentiable in some interval of values of $q$ and moreover it is there strictly convex, we have a large deviation principle like (\ref{lhi}), eventually slightly modified as in  (\ref{secld}). Actually, Remark \ref{rid} becomes particularly pertinent in view of the next result by Hunt and Kaloshin. They in fact showed  
  that for a prevalent set of $C^1$ observables $f:\mathbb{R}^n \to \mathbb{R}^m$, and for $1\le q\le 2$, $D_q^f$ is given by 
 \begin{equation}\label{H33}
 D_q^f=\min(D_q,m).
 \end{equation}
 
 This result implies that when $m$ is smaller than the $D_q$'s, the image measure is not anymore multifractal, in the sense that   all  the $D_q^f$ are equal to $m$, for $1\le q\le 2.$ The function $(q-1)D_q^f$ is not strictly convex and therefore no large deviation principle holds for $\D$.
 
 On the other hand, if  $m$  is larger than the $D_q$, we have from equation (\ref{H33}) that $D_q^f=D_q$, at least for $1\le q\le2$ . Therefore, the image measure inherits some part of the generalized dimensions spectrum from the original measure, which could influence the fluctuations of $\D$ around $D_1$.

 Apart the threshold imposed by $m$,  the observable $f$ will not exhibit itself explicitly in the  detection of  $\D$ given by  equation (\ref{H33}) in the range $q\in[1,2].$ One could ask if the influence of $f$ will manifest for values  of $q$ outside the interval $[1,2]$. Hunt and Kaloshin gave  examples of dynamical systems where $D_q\neq D^f_q$ for $q\notin [1,2].$ \\We will instead  show that in presence of non prevalent observable the image measure will not in general satisfy Eq. (\ref{H33}). In conclusion: {\em the signature of the observable $f$ could become apparent by affecting the typical value of $\D$ for large and small $q$, or when $f$ is not prevalent}\footnote{We remind however that the observable manifests itself in the computation and in the detection of the extremal index, as we showed in formula (\ref{gen}).}. This issue could be important when we analyze time series generated by physical observables, especially if the underlying dynamical systems is high dimensional. We will study a few of those cases in a future publication.

\subsection{Examples}
In the following examples we will mostly consider the baker's map for which we can establish rigorous results. The baker map does not give a conformal repeller, but we could reduce to it by conditioning on the invariant manifolds. 

\subsubsection{Vertical linear observable} We consider the baker map studied in the previous sections. The attractor of this map has a multifractal structure \cite{ott,VVV}.\\

Let us take the observable $f(x,y)=x$ and consider
$$
\int f_*\mu(B(z, r))^{q-1} df_*\mu(z),
$$
where $z\in \mathbb{R}$.\\
First, by definition of image measure we bring the integration over the  SRB measure supported on the  baker's attractor:
$$
\int f_*\mu(B(f(\overline{v}), r))^{q-1} d\mu(\overline{v}), \ \overline{v}=(x,y).
$$

We notice that $f_*\mu(B(f(\overline{v}), r))$ is exactly the SRB measure of  a vertical strip centered at $x$ and with width $r$. We now use disintegration and  write
$$
\Sigma(x):=\{(x,y); |x-y|<r, y\in [0,1]\}=f^{-1}(B(x,r)),
$$
and
$$
\int f_*\mu(B(x, r))^{q-1} d\mu(\overline{v})= \int_{\mathcal{F}_s} \int_{W_{s,\nu}}f_*\mu(B(x, r))^{q-1} d\mu_{s,\nu}d\zeta(\nu),
$$
where $W_{s,\nu}$ denotes an horizontal stable manifold indexed with $\nu$ and $\zeta$ is the counting measure over the stable foliation $\mathcal{F}_s.$ Since stable manifolds are horizontal  segments of length  $1$ emanating from all but countably many points $y$ on the $y$-axis, we will, from now on, identify $\nu$ with $y$ and the first integral on the right hand side of the expression above will be evaluated between $0$ and $1.$ It has been proved in section 4.2 that
$$
f_*\mu(B(x, r))=\mu_{s,y}(\Sigma(x)\cap W_{s,y}),
$$ 
which is the conditional measure of a ball of radius $r$ around the point with abscissa $x$. This measure does not depend on $y,$ and also the conditional measures are the same on all stable fibers and, as we said in section 4.1, they are the invariant ({\em balanced}), measure of a one-dimensional conformal repeller with two linear branches of slopes $\lambda_1^{-1}$ and $\lambda_2^{-1}$ and weights $\alpha$ and $1-\alpha.$
 In conclusion
$$
D_q^f=\lim_{r\to0}\frac{\log \int f_*\mu(B(x,r))^{q-1}d f_*\mu(x)}{(q-1)\log r}= \lim_{r\to0}\frac{\log \int_0^1 \mu_{s,y}(\Sigma(x)\cap W_{s,y})^{q-1}d\mu_{s,y}}{(q-1)\log r}
$$

Therefore the generalized dimensions spectrum $D_q^f$ of the image measure will be that of the associated $1-D$ IFS which are the solution of the transcendental equation \cite{ott, VVV}:
$$\alpha^q\lambda_a^{(1-q)D^f_q}+(1-\alpha)^q\lambda_b^{(1-q)D^f_q}=1,$$

 which differs from the generalized dimensions $D_q$ of the baker attractor (in fact we have $D_q=1+D_q^f$) \cite{ott}. Therefore, this observable does not belong to the prevalent set of the Hunt-Kaloshin theorem, but we already proved that the observable $f(x,y)=x$ is not prevalent.
 \subsubsection{Horizontal and oblique  linear observable}
 
We take now first the observable $f(x,y)=y$ (which is prevalent for $\D$) and we disintegrate along the unstable manifolds $W_{u, \iota}$ (see section 4.1), where the index $\iota$ characterizes the uncountable family of unstable leaves which foliate  baker's attractor. In this case we move up a horizontal strip of width $r$: $\Sigma(y):=\{(x,y); |x-y|<r, x\in [0,1]\}=f^{-1}(B(f(y),r)).$  This strip has a measure which is independent of its height and of the unstable leaf $W_{u, \iota}$; it is therefore given by   $2r$ (the vertical thickness) times $1$ which is the full balanced measure along the $x$-axis.  Remember also that each unstable manifold carries a normalized Lebesgue measure  $\text{Leb}$. Therefore we have
$$
D_q^f=\lim_{r\to0}\frac{\log \int f_*\mu(B(x,r))^{q-1}d f_*\mu(x)}{(q-1)\log r}= \lim_{r\to0}\frac{\log \int \text{Leb}(\Sigma(y)\cap W_{u,\iota})^{q-1}d\text{Leb}(y)}{(q-1)\log r}=
$$
$$
\lim_{r\to0}\frac{\log \int (2r)^{q-1}d\text{Leb}(y)}{(q-1)\log r}=1,
$$
which shows that all the generalized dimensions for the image measure are equal to $1,$ and {\em this proof works for any $q.$}
The same proof immediately generalizes to  linear scalar observables of the form $f(x,y)=ax+by+c$, with $b\ne 0$ and it will give that the $D^f_q=1, \forall q, $ just establishing that there are no deviations from the typical value $D_1=1.$ 

\subsubsection{Numerical verification} 
\begin{figure}[h!]
    \centering
     \includegraphics[height=2.5in]{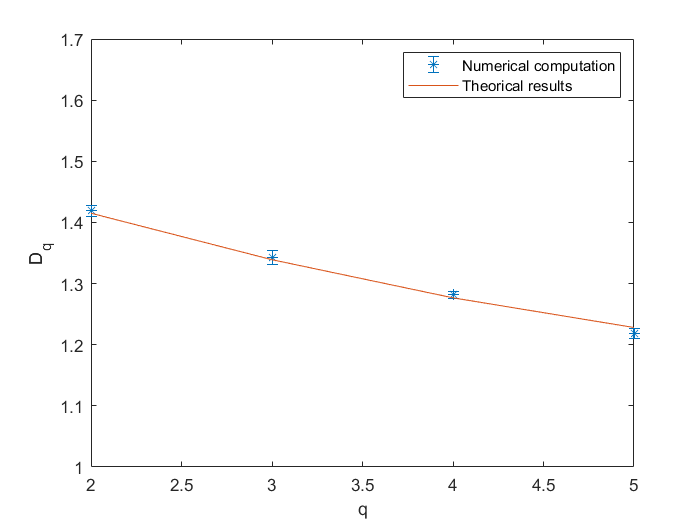}
   
   \caption{$D_q^f$ computed for a Sierpinski gasket and the observable $f(x,y)=(0.2x+2y,x^2),$ compared with theoretical values of $D_q$.}
\label{serp}
\end{figure}

We computed numerically the $D_q^f$ of the baker map for the linear observables introduced at the end of the previous section, using the EVT based method developed  in \cite{dq}. For different values of $a$ and $b\neq 0,$ we found a spectrum of generalized dimensions very close to $1$ up to $q=5$ (the discrepancy of the method for high $q$ yields imprecise results for $q>5$). In conclusion, we believe that for a prevalent set of smooth scalar observable (or more generally for those for which the dimensionality $m$ is smaller than the generalized dimensions of the system), $D_q^f$ will be $1$ (or $m$) for all $q$. We successfully tested this matter numerically for different $C^1$ scalar observables. Similar results are found for the H\'enon system and for a multifractal Sierpinski gasket that we constructed with the iterated functions system technique presented in \cite{dq}. We took the probabilities $p_1=p_2=1/4$ and $p_3=1/2$. The generalized dimensions for this system are explicit (see \cite{dq}) and  comprised between $1$ and $2.$ When  we take an observable in $\mathbb{R}^2$, we find a perfect agreement between $D_q^f$ and $D_q$ (see figure \ref{serp}) for $q$ ranging from $2$ to $5$. This is a sign that the result of Hunt-Kaloshin may hold for a large class of systems in a much broader range for $q$ than the interval $[1,2]$. We proceeded our computations
using the EVT based method, as for the baker map. We took trajectories of length $10^8$
and blocks of size $10^4$.

\section{Applications}\label{clim}

Extreme Value Theory for dynamical systems has been a promising framework for devising metrics to study the climate system. Several heuristic studies~\cite{nature,dq,FF,messori} have focused on the applicability  to climate variables. Although those studies did not have an immediate mathematical justification,   they provided insights on the multifractal, non-stationary nature of the climate attractor. Here  we give an a posteriori justification of those results,  showing that the EVT can be applied to a wide range of observables. It is worth mentioning that in the physical applications, instead of looking at an observable defined on the phase space, we will follow it in time,which would be equivalent by assuming ergodicity of the transformations.
\subsection{From scalar to vector-valued observables}
A common approach to compute the dimensions of the attractor  is to use embedding techniques. It consists in taking a $C^2$ scalar observation $\alpha(x)$, that is accessible through measurement, and constructing a delay observable $f_k(x)=(\alpha(x),\alpha(T(x)),\dots,\alpha(T^{k}(x)))$ (a lag parameter is sometimes added in the numerical studies) \cite{em}. Takens' Theorem states that when the map $T$ is a $C^2$ diffeomorphism defined on a smooth manifold of dimension $D$ and  $\alpha$ is a $C^2$ function, the observable $f_k$ is {\em generically} an embedding into $\mathbb{R}^{2D+1}.$  The notion of (topological) genericity echoes with the notion of prevalence  used in the Hunt-Kaloshin Theorem. In fact  Takens' Theorem has been strengthened in \cite{em} just by using prevalent observables defined directly on a compact invariant subset of some $\mathbb{R}^l$ and where $D$ is now the box counting dimension of that compact set.  As $f_k$ is a diffeormorphism, it preserves the fine structure of the attractor, it allows to reconstruct it and its dimension can be computed numerically. The results of Hunt-Kaloshin, however, state that it is enough that $k>D_1$ to have access to the local dimensions of the attractor, provided the delay coordinate observable $f_k$ is prevalent (which is surely the case for a dense set of $C^1$ observables  $\alpha$). Therefore the Hunt-Kaloshin result provides a more efficient way to access the dimension of the attractor, but it is surely not enough, in general, to reconstruct it. 
\begin{remark}
It is very important to point out that in order to get the dimension we {\bf do not need to reconstruct the attractor}, since the dimension of the image measure is provided directly by the extreme value technique as one of the parameter in the numerical detection of  Gumbel's law. 
\end{remark}

It is well known that embedding techniques often work efficiently for a number of delay coordinates $k$ that is much smaller than the theoretical value of $2D$ prescribed by the Takens theorem, at least when it comes to the computations of the attractor dimensions \cite{em}. For the Lorenz system for example, it is enough that $k=3$ \cite{sauer2}. This result is particularly well understood with the Hunt-Kaloshin results, for which it suffices  to have $k>D_1$ to have that $\D=D_1$ almost everywhere. To illustrate this fact let us begin with a more general consideration. Let  $f$ be a $C^1$ scalar function defined on a neighborhood $\mathcal{U}$ of our attractor and let the map $T$ be smooth enough in order to apply Theorem \ref{HUU}. We then define the vector-valued function, with values in $\mathbb{R}^k$ and components $(f_1, \dots, f_k):$ 
\begin{equation}\label{vbn}
f_j(x)=f(T^{j-1}x), \ x\in \mathcal{U}, \ j=1, \dots k.
\end{equation}
We are now in position to apply the Hunt-Kaloshin Theorem and to look at the least embedding dimension. We tested it on the Lorenz map introduced in section 4.3, reconstructed with the Euler method with step $h=0.01$. With this iterative procedure, the system could be seen as a discrete mapping, to which we apply our theory.\\

\begin{figure}[h!]
    \centering
     \includegraphics[height=2.5in]{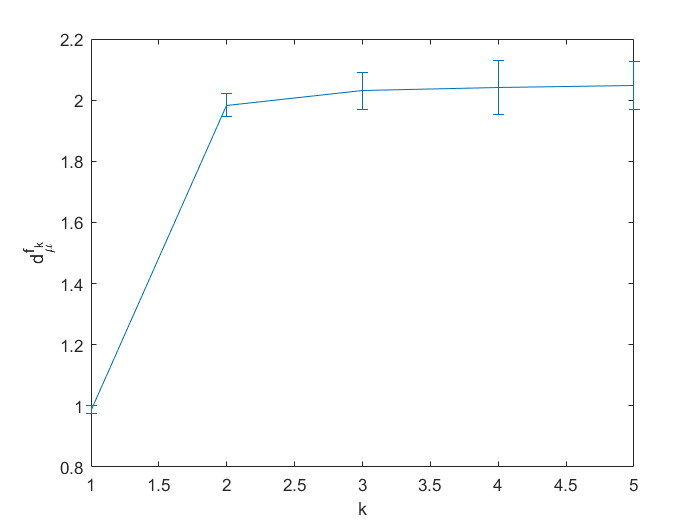}
   
   \caption{Value of $d_{\mu}^{f_k}$ found for different values of $k$, for the Lorenz system and the scalar observable $f(x,y,z)=x$. The parameters used are described in the text. The error bars are the standard deviations of the results.}
   \label{loremb}
\end{figure}

 In figure \ref{loremb}, we computed the local dimensions associated with a vector-valued function like Eq. (\ref{vbn}) with
  $f$  defined by the projection on the $x$ axis, $f(x,y,z)=x$. For the computations, we generated trajectories of $4\cdot 10^7$ points and took a block size of $10^4$. The results are averaged over $20$ different trajectories and target points. We find indeed that $\D=k$ for $k<D_1\approx 2.04$ and becomes constant equal to $D_1$ for $k>D_1$. This suggests that the EVT based methods to compute local dimensions are suited to determine the dimensions of the attractor from a scalar observation $f$: it is enough to construct the delay-coordinate observable $f_k$ and compute its associated dimension $d_{\mu}^{f_k}$ for different values of $k$, until they do not vary anymore as we increase $k$ or until a non-integer value is obtained. We have then attained the dimension of the attractor $D_1$.

\subsection{Vector-valued observables}
At the end of the previous section we showed how to construct a vector-valued observable by composing a given scalar function with the dynamics. This is the spirit of the embedding approach. It turns out that this procedure has its limits, for instance if the map $T$ is not  regular enough to construct prevalent delayed-coordinate observables. 
Sometimes a variety of observables are available and could be used to compute the local dimensions of the original system, provided their cardinality, say $L$ is large enough. Physicists can  measure various quantities associated to the system (temperatures, pressures, velocities, positions\dots), which could be arranged as the outcomes of a function with values in $\mathbb{R}^L.$ It is enough that $L$ is larger than the information dimension of the system to be able to compute the latter. For example, in order to study the dynamics of the atmospheric circulation over the North Atlantic, several authors \cite{nature,dq,FF,messori} have considered an observable that is a vector containing the values of the sea-level atmospheric pressure on a grid of $\approx 10^3$ locations over the North Atlantic. These data were analyzed by \cite{nature} and \cite{dq}.
The source of variability of the sea-level pressure (SLP) atmospheric pressure data has been related to the properties of the atmospheric circulation, namely the switching between different weather regimes~\cite{FF}, the occurrence of extreme weather events~\cite{messori} and the non-stationarity of the underlying attractor due to climate change~\cite{nature}. The local dimensions computed from these observations were centered around the value 13, which is much smaller than the dimension of the space where evolves the observable ($\approx 10^3$) \cite{nature}. The Hunt-Kaloshin theorem can provide a justification to these results: it is enough that the dimension of the ambient space $k$ where  the observable evolves  (here $\approx 10^3$) is larger than $D_1$ (the information dimension of the underlying system) to get that $\D$ is equal to $D_1$ almost everywhere. In other words, when we compute some low and non-integer values for $\D$, it could be the sign that the information dimension of the underlying system is also not integer and much smaller than that  of the ambient space.\\ 

 To study what happens when the dimension of the observable is smaller than the information dimension of the attractor, we use again the SLP  data from the reanalysis of the National Centers for Environmental Prediction (NCEP) \cite{kalnay}, but we  now investigate the effect of averaging the information over all the grid points. To this purpose, we define the function
$$\phi_Z(X)=-\log |\langle SLP(Z)\rangle - \langle SLP(X)\rangle|,$$

where $\langle SLP(X)\rangle$ is the spatial arithmetic average value of the $x_i$:
$$\langle SLP(X)\rangle=\frac1n\sum_{i=1}^nx_i,$$

and $Z=(z_1,\dots,z_m)$ corresponds to a particular configuration of the pressure field.

\begin{figure}[h!]
    \centering
    \begin{subfigure}[t]{0.5\textwidth}
        \centering
\includegraphics[height=2.5in]{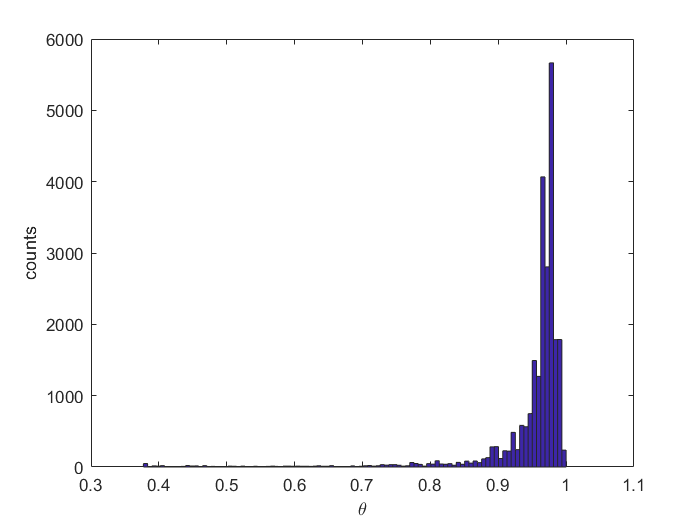}
\label{thetadistrib}
    \end{subfigure}%
    ~
    \begin{subfigure}[t]{0.5\textwidth}
        \centering
     \includegraphics[height=2.5in]{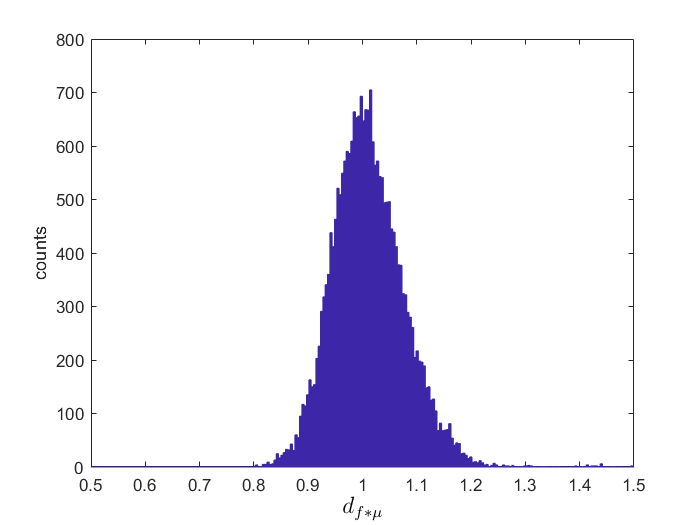}
     \label{d}
    \end{subfigure}
   \caption{Distributions for the atmospheric circulation (SLP) data presented in the text of the EI (left) and of the $d_{f_*\mu}(f(z))$ (right) found for different points $z$ of the attractor.}
   \label{fig4}
\end{figure}

To compute $\theta$ and $\D$ associated with this observable, we perform a computation of the empirical distribution of the variable $M_{50}$ defined in  Eq. (\ref{MAX}), for different points $Z\in X$. For each of them, we find that the best fit of the empirical distribution is a Gumbel law of scale parameter close to $1$ and estimates of the extremal index are close to $1$, like in the baker map situation. The fittings are performed using the Matlab function {\em gevfit}. For the computation of the extremal index, we used the estimate $\hat\theta_0$ introduced in \cite{ei}, with a threshold value equal to the 0.99-quantile of the observable distribution. Both the values of the extremal index and of $d_{f*\mu}$ that we found  have small variability around $1$, due to finite effects. In figure \ref{fig4} the distributions of the values found for $\theta$ and $d_{f*\mu}$ over the different points $Z$ of the attractor are represented. These estimates are in perfect agreement with the results presented in this work.

%

We now study the statistics of the number of visits of the observable in the neighborhood of a particular value for the presented climate data, as we explained  in section \ref{viss}. We can observe in figure \ref{fig4} on the left, that finite effects lead to an estimate of $\theta$ that is slightly smaller than $1$, due, among other reasons, to the persistence of the orbits in the neighborhood of the point $z$. This clustering is very likely to disappear if the amount of data  allows to take a higher threshold, and we would observe a pure Poisson distribution at the limit of high threshold. We get a distribution that is very close to P\`olya-Aeppli of parameters $t$ and the extremal index computed at a finite resolution. This is consistent with the discussion on the figure \ref{comp1}: in many situations, although we know we do not have a P\`olya-Aeppli distribution, it seems that it still models the limit law quite well. In our computations, we studied visits of the observable $f(X)=\frac1n\sum_{i=1}^n x_i$ in the ball $(f_0-r,f_0+r)$, where $f_0$ is the value taken by the observable on July, $1^{st}$, 1948 and $r=e^{-u}$, $u$ being the $0.98-$ quantile of the distribution of $\phi$.\\

\begin{figure}[h!]
    \centering
     \includegraphics[height=2.5in]{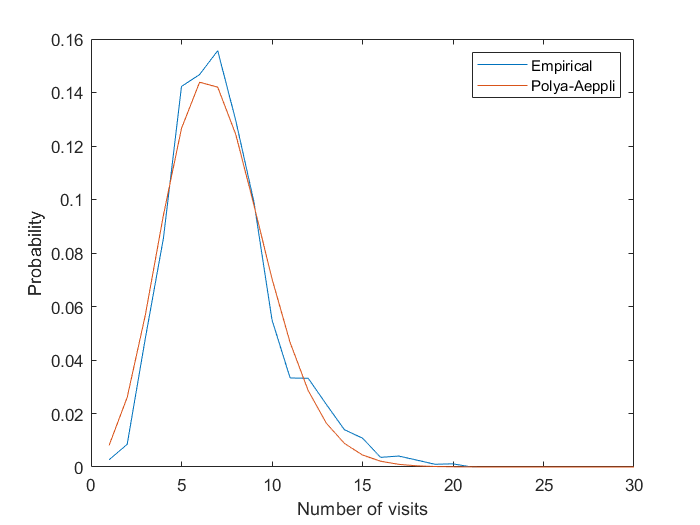}
   
   \caption{Comparison between the empirical distributions of the number of visits of the observable mean value in a ball centered at $f(z)$ and a Polya-Aeppli distribution for the climate data presented in the text.}
   \label{fig5}
\end{figure}

The existence of extreme value laws for recurrences of physical observables also justifies the results obtained in~\cite{farandavaientigrl,farandaetalclimdyn,alvarez,gualandi,rodrigues}, where the rate of statistical convergence to the extreme value laws was used to estimate the characteristic recurrence time of temperature values (termed recurrence spectra in~\cite{alvarez}). They show that despite the slow convergence of the dynamical systems metrics towards unknown asymptotic values, their distribution is reminiscent of an underlying  high-dimensional  attractor. On this object,  the recurrences around high dimensional fixed or periodic points determine interesting dynamical behaviors such as switching between metastable states~\cite{nature}, critical phenomena~\cite{gualandi} or different basin of attraction~\cite{brunetti} that can be detected by deviations of the dynamical indicators from their expected asymptotic behavior. 
\section{Conclusions}
 This paper contains a few rigorous results illustrated by several examples. The latter  are worked out  relatively easily, but it was important for us to show that the statistical indicators established by the theory can be explicitly computed and compared with the numerical simulations. The dynamical systems we considered have strong mixing properties, in particular they exhibit exponential decay of correlations on suitable spaces of observables. This  allowed us to use a very efficient perturbative theory and  compute the extremal index in a broad variety of situations \cite{ei}. We believe that our results  could be generalised to larger class of systems, even non-uniformly hyperbolic, or exhibiting intermittency \cite{PDC1, PDC2},  using for instance techniques with more probabilistic flavor \cite{FFT, FFT13, Ay, freitascat, book}.  In this perspective,   we also considered more complex systems and physical time series that we analyzed numerically and that  can be tested and interpreted in the framework of our theory. As we explicitly shown in the last section,  we believe that our results are  useful for physical and natural systems, in the sense that they provide a formal framework for the applications presented in~\cite{nature,nature2,brunetti,rodrigues,gualandi}. They partially answer the concern raised in~\cite{buchow} about the slow convergence of dynamical system metrics for climate data and make useful asymptotic theorems for finite data sets.
\section{Acknowledgements}
S.V. thanks the Laboratoire International Associ\'e LIA LYSM, the
INdAM (Italy), the UMI-CNRS 3483, Laboratoire Fibonacci (Pisa) where this work has
been completed under a CNRS delegation, the {\em Centro de Giorgi} in Pisa and  UniCredit Bank R\&D group for financial support through the
Dynamics and Information Theory Institute at the Scuola Normale Superiore.
 SV thanks N. Haydn for enlightening discussions on this paper.
The authors thank A. Coutinho and W. Bahsoun for useful exchanges on a few part of this article.


\end{document}